\RequirePackage{fix-cm}
\documentclass[smallextended]{svjour3}
\usepackage[utf8]{inputenc}
\usepackage[T1]{fontenc}
\usepackage[english]{babel}
\usepackage[breaklinks,colorlinks=true,bookmarksnumbered,
bookmarksopen=true,citecolor=red,urlcolor=blue]{hyperref}
\smartqed  % flush right qed marks, e.g. at end of proof
\usepackage{enumerate}
\usepackage{amsmath,amssymb,amsfonts,latexsym,mathptmx,mathtools}
\usepackage{bbm}
\usepackage{url}
\usepackage{datatool,filecontents}
\usepackage{pgfplots}
\usepackage{tikz}
\usetikzlibrary{patterns}
\usepackage{varwidth}
\usepackage{xcolor}
\usepackage{graphicx}
\usepackage{url}
\newcommand{\R}{\mathbb{R}}
\newcommand{\N}{\mathbb{N}}
\newcommand{\B}{\mathbb{B}}
\newcommand{\D}{\mathcal{D}}
\newcommand{\U}{\mathcal{U}}
\newcommand{\A}{\mathcal{A}}
\newcommand{\K}{\mathcal{K}}
\newcommand{\Z}{\mathcal{Z}}
\newcommand{\ZN}{\mathcal{Z}_0}

\newcommand{\di}{\textnormal{dist}}
\newcommand{\dom}{\textnormal{dom}}
\newcommand{\gph}{\textnormal{gph}}
\newcommand{\rge}{\textnormal{rge}}
\DeclareMathOperator*{\Limsup}{Limsup}

\newcommand{\cl}{\mathop{\textnormal{cl}}}
\newcommand{\inter}{\mathop{\textnormal{int}}}

\newcommand{\bd}{\mathop{\textnormal{bd}}}
\newcommand{\Reg}{\textnormal{Reg}}
\newcommand{\cone}{\textnormal{cone}}
\journalname{}

\begin{document}
\title{The Least Singular Value Function in Variational Analysis\thanks{Research of the second author was partly supported by the US National Science Foundation under grant DMS-2204519 and by the Australian Research Council under Discovery Project DP250101112.}}

\author{Mario Jelitte \and Boris S. Mordukhovich}

\institute{Mario Jelitte \at 
Faculty of Mathematics, Technische Universit\"{a}t Dresden, 01062 Dresden, Germany\\
\email{Mario.Jelitte@tu-dresden.de}
\and
Boris S. Mordukhovich \at 
Department of Mathematics, Wayne State University, Detroit, Michigan 48202, USA\\
\email{aa1086@wayne.edu}}

\maketitle

\begin{abstract}
Metric regularity is among the central concepts of nonlinear and variational analysis, constrained optimization, and their numerous applications. However, metric regularity can be elusive for some important ill-posed classes of problems including polynomial equations, parametric variational systems, smooth reformulations of complementarity systems with degenerate solutions, etc. The study of stability issues for such problems can often not rely on the machinery of first-order variational analysis, and so higher-order regularity concepts have been proposed in recent years. In this paper, we investigate some notions of mixed-order regularity by using advanced tools of first-order and second-order variational analysis and generalized differentiation of both primal and dual types. Efficient characterizations of such mixed-order regularity concepts are established by employing a fresh notion of the least singular value function. The obtained conditions are applied to deriving constructive criteria for mixed-order regularity in coupled constraint and variational systems.\vspace*{-0.05in}

\vspace*{-0.05in}
	
\end{abstract}\vspace*{-0.2in}

\section{Introduction}\label{sec.Intro}\vspace*{-0.1in}

It has been well recognized that \textit{metric regularity} of set-valued mappings is among the fundamental notions of variational analysis and related topics. In particular, this concept has been used
to establish calculus rules of generalized differentiation, to formulate verifiable necessary optimality conditions for constrained optimization problems, to tackle Lipschitzian stability and sensitivity issues for generalized equations and equilibrium problems, variational and quasi-variational inequalities, etc.; see, e.g., \cite{BSh:PAOP:00,Ioffe,Mo06,Mo18,Mo24,Rocka,thibault} and the references therein. We specially highlight that metric regularity and the associated notions of Lipschitzian stability play a crucial role in the design and justification of Newton-like algorithms to solve subgradient inclusions and other classes of optimization-related problems, which is demonstrated in \cite{BSh:PAOP:00,Rock_ImpF,FP03-1,ISo:NTMOVP:14,KMP23,KK02,Mo24} among many publications in this direction.

On the other hand, metric regularity may not hold, or be largely restrictive, in numerous settings important for variational analysis and optimization. Let us mention, in particular, the failure of metric regularity for solution mappings associated with parametric variational inequalities and complementarity systems; see \cite{Mo08} and \cite[Section~3.3]{Mo18} for more details. Having this in mind, we focus in what follows on the study of less restricting regularity notions, which provides anyway a sufficient amount of information for variational theory and applications beyond the realm of metric regularity. The notions of our study can be generally called ``mixed-order regularity" that extend and unify the two groups of relaxed regularity conditions. The first one was introduced in \cite{Av85,Tret84} under the name of 2-{\em regularity} and then had been successfully employed in developing  numerical methods of optimization \cite{Iz_CrSol_16,FIJ23B}, stability and sensitivity of local minimizers \cite{AruIz20,FIJ21,GM,Iz_CrSol_16,Mo24}, error bounds and 
higher-order necessary optimality conditions  \cite{AAI05,GM,J24,Tret84}, etc.

The second newer notion of mixed-order regularity, introduced in \cite{Gf13} and studied in our paper, is Gfrerer's {\em metric pseudo-regularity}; see, e.g., \cite{BM25,GfOu16,JDiss24,J24} for more recent developments based on methods of second-order variational analysis.

To bring the above regularity concepts under the same roof, the first author has recently introduced in \cite{JDiss24} the \textit{least singular value} (LSV) function, which is investigated and applied in the current paper via first-order and second-order primal and dual constructions of generalized differentiation. We show that this approach allows us, from one side, to characterize the classical notions of metric regularity and related properties, while from the other side, it leads us to novel characterizations of the mixed-order regularity notions mentioned above. 

The rest of the paper is organized as follows. Section~\ref{sec.2} recalls some basic definitions from variational analysis and then formulates and discusses the major construction of the least singular value function and the associated notion of {\em singularity}. In Section~\ref{sec.2.2}, we establish several 
estimates for {\em primal-type} derivatives of the LSV function.

Section~\ref{sec.MR} presents some review of dual-type constructions of generalized differentiation revolving around {\em coderivatives} of set-valued mappings which are instrumental to derive complete characterizations of metric regularity and related well-posedness properties. Now we look at these issues from the viewpoint of the LSV function, which provides a bridge to study mixed-order regularity. 

Section~\ref{sec.SC} is devoted to the notion of {\em metric 2-regularity} and plays the central role in the paper being based on the LSV function. We derive here efficient characterizations of metric 2-regularity in terms of the graphical derivative and mainly via the coderivative of multifunctions while providing also various specifications and examples. The subsequent Section~\ref{sec.Gfr} establishes precise relationships between metric 2-regularity and Gfrerer's regularity concept. This leads us to deriving characterizations as well as sufficient conditions for the latter from those established in the preceding section.

The next two sections provide applications of the obtained results for mixed-order regularity to two structured classes of systems that overwhelmingly appear in problems of variational analysis and optimization. The first class of {\em parametric constraint systems} is studied in Section~\ref{sec.App}, and the second class of {\em variational systems} is considered in Section~\ref{sec.App.VS}. The results established in these two sections are most based on the coderivative construction, which enjoys comprehensive calculus rules. 

The concluding Section~\ref{sec:conc} summarizes the main achievements of the paper and discusses some problems of the future research.\vspace*{-0.25in}

\section{Basic Definitions and Preliminaries}\label{sec.2}\vspace*{-0.1in}

In this section, we first present some basic notions of variational analysis broadly used below and then proceed with considering generalized least singular values of set-valued mappings that depend on two variables being positively homogeneous in one of them. The presented material serves as a preparation for our subsequent analysis, where least singular values are employed to characterize certain regularity concepts.\vspace*{0.02in} 

Recall some definitions from \cite{Rocka} that are used in what follows. Throughout the paper, we work in finite-dimensional  Euclidean spaces, where $\Vert\cdot\Vert$ stands for the Euclidean norm, and where $\di[\xi,Z]$ is the distance of a point $\xi$ to a set $Z$, with the understanding that the distance to the empty set is $\infty$. The closed unit ball is $\B$, and the unit sphere is $\bd(\B)$, where $\bd(Z)$ stands for the boundary of $Z$. By $\inter(Z)$ and $\cl(Z)$, we denote interior and closure of $Z$, respectively. The conic hull of $Z$ is
\begin{align*}
\cone(Z):=\bigcup_{\gamma\geq 0}\gamma\cdot Z=\bigcup_{\gamma\geq 0}\{ \gamma\cdot z\mid z\in Z\},
\end{align*}
and the (Bouligand-Severi) \textit{tangent/contingent cone} associated with $Z$ at $\xi\in Z$ is
\begin{align}\label{tan}
T_Z(\xi):=\Limsup_{\tau\searrow 0}\frac{Z-\xi}{\tau},
\end{align}
where the outer limit, $\Limsup$, is understood in the sense of Painlev\'{e}-Kuratowski.

For a set-valued mapping $H:\R^p\rightrightarrows\R^q$, we denote its graph, domain, and range, by $\gph\,H$, $\dom\,H$, and $\rge\,H$, respectively. $H$ is \textit{positively homogeneous} if
\begin{align*}
H(\gamma\cdot z)=\gamma\cdot H(z)=\{\gamma\cdot \eta\mid \eta\in H(z)\}\;\mbox{ for all }\;\gamma>0,\;z\in\R^p,
\end{align*}
and it is \textit{outer semicontinuous} at a point $z\in\dom\,H$ relative to a set $Z$ if
\begin{align*}
\Limsup_{z'\overset{Z}{\rightarrow} z} H(z')\subset H(z).
\end{align*}
If the latter holds for each $z\in\R^p$ relative to $\R^p$, then $H$ is simply called outer semicontinuous. It is easy to see that the outer semicontinuity of a set-valued mapping is equivalent to the closedness of its graph.

Given a single-valued mapping $\mathcal{A}\colon\R^n\to\R^p$, the \textit{semiderivative} of $\mathcal{A}$ at 
$\xi\in\inter(\dom\,\mathcal{A})$ for $\omega$ is defined by
\begin{align*}
\mathcal{A}'(\xi;\omega):=\lim_{\substack{\tau\searrow 0 \\ \omega'\rightarrow \omega}}\frac{\mathcal{A}(\xi+\tau\omega')-\mathcal{A}(\xi)}{\tau}.
\end{align*} 
If the semiderivative exists, then $\mathcal{A}$ is called  \textit{semidifferentiable at $\xi$ for $\omega$}. As well known, the semidifferentiability of $\mathcal{A}$ at $\xi$ for $\omega$ is equivalent to the existence of a continuous, positively homogeneous mapping $B$ satisfying $B(\omega)=\mathcal{A}'(\xi;\omega)$ and
\begin{align}\label{f.expansion}
\Vert \mathcal{A}(\xi+\tau\omega')-\mathcal{A}(\xi)-\tau B(\omega')\Vert=o(\tau\Vert\omega'\Vert)\quad\textnormal{as }\tau\searrow 0,\,\omega'\rightarrow\omega.
\end{align} 
Observe also that if $\mathcal{A}$ is (Frech\'{e}t) differentiable at $\xi$, then $\mathcal{A}'(\xi;\omega)=\mathcal{A}'(\xi)\omega$ for all $\omega$, where $\mathcal{A}'(\xi)$ is the Jacobian of $\mathcal{A}$ at $\xi$. 

Having further an extended-real-valued function $\varphi:\R^n\to\overline{\R}:=\R\cup\{\infty\}$, the \textit{subderivative} of $\varphi$ at $\xi\in\dom\,\varphi$ in the direction $\omega$ is
\begin{align*}
d\varphi(\xi)(\omega):= \liminf_{\substack{\tau\searrow 0 \\ \omega'\rightarrow \omega}}\frac{\varphi(\xi+\tau\omega')-\varphi(\xi)}{\tau}.
\end{align*}
The following major definition is taken from  \cite{JDiss24}.\vspace*{-0.05in} 

\begin{definition}\label{def.LSV} Given a set $\D\subset\R^n\times\R^p$, consider a mapping $\varXi:\D\times\R^m\rightrightarrows\R^q$ such that $\varXi(\xi,\cdot)$ is positively homogeneous for all $\xi\in\D$. Then the function $\ell_\varXi:\R^n\times\R^p\rightarrow[0,\infty]$ defined by
\begin{align}\label{LSV}
\ell_\varXi(\xi):=\left\{ 
\begin{array}{c l}
\infty & \quad\textnormal{ if }\xi\notin \D,\\
\inf\limits_{\Vert z\Vert=1} \di[0,\varXi(\xi,z)] & \quad\textnormal{ if }\xi\in \D\\
\end{array}\right.
\end{align}
is called the \textit{least singular value $($LSV$)$ function} of $\varXi$. A point $\bar{\xi}\in \D$ is called \textit{$\varXi$-singular} if the inclusion $0\in\varXi(\bar{\xi},z)$ is satisfied for some $z\neq 0$.
\end{definition}\vspace*{-0.05in}

The name ``least singular value function" is motivated by the special case $\varXi(\xi,z)=\A(\xi)^\top z$ for a matrix-valued mapping $\A:\D\rightarrow\R^{m\times q}$, where the representation
\begin{align*}
\ell_\varXi(\xi)=\min_{\Vert z\Vert=1}\left\Vert \A(\xi)^\top z\right\Vert
\end{align*}
holds for any $\xi\in \D$. Here the LSV function $\ell_\varXi(\xi)$ corresponds to the least singular value of the rectangular matrix $\A(\xi)$, i.e., the least eigenvalue of the positive-semidefinite square matrix $\A(\xi)\A(\xi)^\top$.

As we see below, the LSV function is fundamental for characterizing generalized regularity concepts for mappings when $\varXi$ is associated with some generalized derivative. It follows from \cite[Proposition~4A.6]{Rock_ImpF} that
\begin{align}\label{f.reciprocal}
|\varXi_\xi^{-1}|^+ 
=\left\{ 
\begin{array}{c l}
\infty & \textnormal{ if }\ell_\varXi(\xi)=0,\\
1/\ell_\varXi(\xi) & \textnormal{ if }\ell_\varXi(\xi)>0
\end{array}
\right.
\end{align}
for the \textit{outer norm} 
\begin{align*}
|\varXi_\xi^{-1}|^+:=\sup_{y\in\B}\sup_{z\in \varXi_\xi^{-1}(y)}\Vert z\Vert
\end{align*}
of the positively homogeneous mapping $\varXi_\xi^{-1}(y):=\{ z\mid y\in\varXi(\xi,z)\}$.
Structured mappings represented as sums of single-valued and set-valued components appear frequently in different frameworks of variational analysis and generalized differentiation. The following result taken from \cite[Theorem~4.2]{JDiss24} deals with the LSV function of such mappings and the corresponding notion of $\varXi$-singularity.\vspace*{-0.07in} 

\begin{theorem}\label{theorem.singular.constraint.mapping} Let 
$\D$ be a nonempty subset of an open set $\U\subset\R^n\times\R^p$, let $\Gamma:\D\times\R^m\rightrightarrows\R^q$ be such that 
$\varGamma(\xi,\cdot)$ is positively homogeneous for all $\xi\in\D$, and let $\varXi:\D\times\R^m\rightrightarrows\R^q$ be given as
\begin{align}\label{f.map.Xi}
\varXi(\xi,z):=\A(\xi) z+\Gamma(\xi,z),
\end{align}
where $\A:\U\rightarrow\R^{q\times m}$. Then $\varXi(\xi,\cdot)$ is positively homogeneous for all $\xi\in\D$ and
\begin{align}\label{theorem.singular.f1}
\ell_\varXi(\bar{\xi})=\inf_{z\in\Z}\di\left[ 0,\A(\bar{\xi})z+H(z)\right]\;\mbox{ for any }\;\bar{\xi}\in \D,
\end{align}
where $H:=\Gamma(\bar{\xi},\cdot)$ and $\Z:=\bd\B\cap\dom\,H$. If furthermore $H$ is outer semicontinuous, then $\bar{\xi}$ is $\varXi$-singular if and only if $\ell_{\varXi}(\bar{\xi})=0$. In this case,
the set
\begin{equation}\label{f.set.Z0}
\ZN:=\left\{ z\in \Z \left| \,0\in \A(\bar{\xi}) z+H(z)\right.\right\}
\end{equation}
is nonempty and closed, and the inclusion $\dom \left(d\ell_\varXi(\bar{\xi})\right)\subset T_\D(\bar{\xi})$ holds.
\end{theorem}\vspace*{-0.3in}

\section{Subderivatives of the LSV Function}\label{sec.2.2}\vspace*{-0.1in}

In this section, we focus on lower estimates for the subderivative of the LSV function. The obtained results play a significant role in the subsequent study of (non-) $\varXi$-singularity. The following lemma proved in \cite{JDiss24} is needed below.\vspace*{-0.05in}

\begin{lemma}\label{lemma.closed.domain.range} Let $H:\R^n\rightrightarrows\R^m$ be an outer semicontinuous and positively homogeneous mapping. Then we have:\vspace*{-0.05in}
\begin{enumerate}
\item[\bf(a)] $\dom\,H$ and $\rge\,H=\dom\,H^{-1}$ are cones. 
\item[\bf(b)] $\gph\,H$ and $H(0)$ are closed cones.
\item[\bf(c)] If $H^{-1}(0)=\{ 0\}$, then $\rge\,H$ is closed. 
\item[\bf(d)] If $H(0)=\{ 0\}$, then $\dom\,H$ is closed.
\item[\bf(e)] It holds that
\begin{align}\label{f.poshom.re,resentation}
\rge\,H=\cone\left(H(0)\cup\left(\bigcup_{\substack{\Vert z\Vert=1,\;z\in\dom\,H}}H(z)\right)\right).
\end{align}
\end{enumerate}
\end{lemma}

Now we obtain two theorems that provide lower estimates for the subderivative of the LSV function under different assumptions. As shown below, these assumptions are generally independent. The first result improves \cite[Theorem~4.4]{JDiss24}. At some stage, their proofs are similar, and so we only present parts of the proof that are new.\vspace*{-0.05in}

\begin{theorem}\label{theorem.subderivative.singular.I}
In the setting of Theorem~{\rm\ref{theorem.singular.constraint.mapping}}, let $\bar{\xi}\in \D$ satisfy the following conditions:
\begin{enumerate}
\item[\bf(i)] $\bar{\xi}$ is $\varXi$-singular.

\item[\bf(ii)] $\Gamma$ is outer semicontinuous.

\item[\bf(iii)] $\A$ is semidifferentiable at $\bar{\xi}$ for $\omega\in\R^n\times\R^p$.

\item[\bf(iv)] The set $\A(\bar{\xi})\dom\,H$ is closed.

\item[\bf(v)] There exists $c\geq 0$ such that for any $z\in \ZN$ we find $\varepsilon,\delta>0$ with 
\begin{align}\label{example.subderivative.approximate}
\varGamma(\xi,z')\cap\delta\B\subset\varGamma(\bar{\xi},z')+c\Vert\xi-\bar{\xi}\Vert\B
\end{align}
whenever $(\xi,z')\in \left(\D\times \bd\B\right)\cap\left((\bar{\xi},z)+\varepsilon\B\right)$.

\item[\bf(vi)] It holds that
\begin{align}\label{example.subderivative.singular.separation}
\left(\A(\bar{\xi}) \dom\,H\right)\bigcap \left(-\varTheta\right)\subset\{ 0\},
\end{align}
where the closed cone $\varTheta$ is defined by
\begin{align}\label{f.subderivative.Theta}
\varTheta:=\cl\left(\cone\left(
\bigcup_{z\in \Z} H(z)\right)\right).
\end{align}	
\end{enumerate}
Then, with the number $c\geq 0$ from ${\rm(v)}$, we have the estimate 
\begin{align*}
d\ell_\varXi(\bar{\xi})(\omega)+c\Vert\omega\Vert\geq \min_{z\in \ZN}\left(\di\left[ 0,\A'(\bar{\xi};\omega) z+\varTheta+
\A(\bar{\xi})\dom\,H\right]\right).
\end{align*}
\end{theorem}
{\bf Proof}. There is nothing to prove when $\omega\notin\dom(d\ell_{\varXi}(\bar{\xi}))$; so suppose that $\omega\in \dom(d\ell_{\varXi}(\bar{\xi}))$. The proof is divided into three main steps, and we only present a detailed verification for the first two steps that are new.\\[0.5ex]
{\bf Step~1}. Let us show that $\ZN$ is nonempty.
According to $\omega\in\dom \left(d\ell_\varXi(\bar{\xi})\right)$, we deduce from  (i) and (ii) that there are sequences $t_k\searrow 0$ and $\omega^k\rightarrow\omega$ with
\begin{align}\label{f.subderivative.singular.f0}
\bar{\xi}+t_k\omega^k\in \D\;\mbox{for all }\;k\in\N\;\mbox{ and }\; t_k^{-1}\ell_\varXi (\bar{\xi}+t_k\omega^k)\rightarrow d\ell_\varXi(\bar{\xi})(\omega)<\infty.
\end{align}
It follows from (ii) that the values of $\Gamma$ are closed. Therefore, there exist sequences $\{z^k\}\subset\bd\B$ and $\{ \eta^k\}$ with $\eta^k\in\Gamma(\bar{\xi}+t_k\omega^k,z^k)$ as $k\in\N$ such that
\begin{align}\label{f.subderivative.singular.II.1}
\begin{array}{r l}
t_k^{-1}\ell_\varXi(\bar{\xi}+t_k\omega^k)=t_k^{-1}\left\Vert \A(\bar{\xi}+t_k\omega^k) z^k+\eta^k\right\Vert\\
=t_k^{-1} \left\Vert \A(\bar{\xi}) z^k+t_kB(\omega^k) z^k+\eta^k\right\Vert+t_k^{-1}o(t_k\Vert \omega^k\Vert)
\rightarrow \gamma
\end{array}
\end{align}
as $k\to\infty$ with the number $\gamma$ defined by
\begin{align}\label{f.finite.subderivative}
\gamma:=d\ell_{\varXi}(\bar{\xi})(\omega)\in[0,\infty),
\end{align}
where the mapping $B:\R^n\times\R^p\rightarrow\R^{q\times m}$ is continuous, positively homogeneous, and  satisfies the conditions $B(\omega)=\A'(\bar{\xi};\omega)$ and \eqref{f.expansion}. Multiplying all the terms in \eqref{f.subderivative.singular.II.1} by $t_k$ and utilizing the continuity of $B$ together with $\Vert z^k\Vert=1$ tell us that
\begin{align}\label{f.subderivative.singular.II.2}
\A(\bar{\xi}) z^k+\eta^k\rightarrow 0.
\end{align}
Without loss of generality, assume that $z^k\rightarrow\widehat{z}\in\bd\B$, and so \eqref{f.subderivative.singular.II.2} yields $\eta^k\rightarrow -\A(\bar{\xi})\widehat{z}$. By (ii), we get $-\A(\bar{\xi})\widehat{z}\in \varGamma(\bar{\xi},\widehat{z})=H(\widehat{z})$, i.e.,
\begin{align}\label{f.subderivative.singular.II.3}
\widehat{z}\in \ZN.
\end{align}
Thus it follows from (vi) that $\mathcal{A}(\bar{\xi})\widehat{z}=0$.\\[0.5ex]
{\bf Step~2}. Here we obtain a lower estimate of $\gamma=d\ell_{\varXi}(\bar{\xi})(\omega)$. Observe from the constructions in Step~1, the convergence of $\{ z^k\}$ to $\widehat{z}$ along with \eqref{f.subderivative.singular.II.3}, $\mathcal{A}(\bar{\xi})\widehat{z}=0$, and (v) that for some $c\geq 0$ and $\delta>0$ we get the inclusion
\begin{align}\label{f.p1.1}
\varnothing\neq \varGamma(\bar{\xi}+t_k\omega^k,z^k)\cap\delta\B\subset\varGamma(\bar{\xi},z^k)+ct_k\Vert\omega^k\Vert\B =H(z^k)+ct_k\Vert\omega^k\Vert\B,
\end{align}
which ensures that $z^k\in \Z$ for all $k$. Furthermore, it follows from Step~1 and \eqref{f.p1.1} that
\begin{align*}
\ell_{\varXi}(\bar{\xi}+t_k\omega^k) & =\di\left[ -\left(\A(\bar{\xi}) z^k+t_kB(\omega^k) z^k\right),\Gamma\left(\bar{\xi}+t_k\omega^k,z^k\right)\right]+o\left( t_k\Vert\omega^k\Vert\right)\\
~ & \geq \di\left[ -\left(\A(\bar{\xi}) z^k+t_kB(\omega^k) z^k\right),H(z^k)\right]-ct_k\Vert\omega^k\Vert+o\left( t_k\Vert\omega^k\Vert\right)\\
~ & \geq h(t_k\omega^k)-ct_k\Vert\omega^k\Vert+o\left( t_k\Vert\omega^k\Vert\right),
\end{align*}
where the function $h$ in the third line is defined by
\begin{align*}
h(\nu):=\inf_{z\in \Z}\left(\di\left[-\left(\A(\bar{\xi})+B(\nu)\right) z,H(z)\right]\right).
\end{align*}
Due to \eqref{f.subderivative.singular.II.3}, we have $\ZN\neq\varnothing$, which yields $h(0)=0=\ell_{\varXi}(\bar{\xi})$. Moreover, $h(\nu)\geq 0$ for any other $\nu$, and so $dh(0)(\nu)\geq 0$ for each $\nu$. Thus it follows from \eqref{f.subderivative.singular.II.1}, \eqref{f.finite.subderivative}, and the estimates above that
$\infty>d\ell_\varXi(\bar{\xi})(\omega)+c\Vert\omega\Vert\geq dh(0)(\omega)\geq 0$.\\[0.5ex]
{\bf Step~3}. It remains to estimate $dh(0)(\omega)$ from below, which was already done in the proof of \cite[Theorem~4.4]{JDiss24}.\qed\vspace*{0.05in}

Observe that condition (i) yields $z\in \dom\,H\setminus\{ 0\}$ with 
$-\A(\bar{\xi})z\in\varTheta$. It follows from (vi) that the latter can only be true when $z\in\ker \A(\bar{\xi})$, which means that the theorem is designed for rather special $\varXi$-singular points.\vspace*{0.05in} 

The conclusion of the theorem can be relaxed as in the following corollary whose simple proof merely combines Theorem~\ref{theorem.subderivative.singular.I} and Lemma~\ref{lemma.closed.domain.range}.\vspace*{-0.05in}

\begin{corollary}\label{c.suff.singular.I}
In the setting of Theorem~{\rm\ref{theorem.singular.constraint.mapping}}, let $\bar{\xi}\in\mathcal{D}$, and let $\omega\in\R^n\times\R^p$ satisfy  conditions {\rm(i)--(vi)}. Then, with the number $c\geq 0$ from {\rm(v)}, we get the estimate
\begin{align*}
d\ell_\varXi(\bar{\xi})(\omega)+c\Vert\omega\Vert\geq \min_{z\in \ZN}\left(\di\left[ 0,\A'(\bar{\xi};\omega) z+\cl(\rge\,H)+
\A(\bar{\xi})\dom\,H\right]\right).
\end{align*}
In particular, the following is sufficient for the fulfillment of {\rm(vi)} in Theorem~{\rm\ref{theorem.subderivative.singular.I}}:
\begin{align}\label{f.suff.cond6}
\left(\mathcal{A}(\bar{\xi})\dom\,H\right)\bigcap\cl(\rge\,H)\subset\{ 0\}.
\end{align}
\end{corollary}

Condition (v) of Theorem~\ref{theorem.subderivative.singular.I} can be interpreted as a variant of the \textit{calmness} property \cite{Rock_ImpF} for the mapping $\varGamma$ with respect to the first variable (uniformly in the second one), while the size of the number $c\geq 0$ therein has a crucial influence on the usefulness of the theorem. The conclusion of the theorem is vacuous when
\begin{align*}
\min_{z\in \ZN}\left(\di\left[ 0,\A'(\bar{\xi};\omega) z+\varTheta+
\A(\bar{\xi})\dom\,H\right]\right)\leq c\Vert\omega\Vert
\end{align*}
since this only leads us to the trivial estimate $d\ell_{\varXi}(\bar{\xi})(\omega)\geq 0$, which is always fulfilled for any $\varXi$-singular point $\bar{\xi}$.
In applications, we will be particularly interested in the circumstances where 
$d\ell_{\varXi}(\bar{\xi})(\omega)>0$. To guarantee this, it is convenient to consider situations with $c=0$.

To check in practice the fulfillment of conditions (ii), (iv), and (v) of Theorem~\ref{theorem.subderivative.singular.I} is challenging. However, this task is simplified when the mapping $\varGamma$ is \textit{polyhedral} in the sense of \cite{Rob1981}, i.e., $\gph\,\varGamma$ is a finite union of (convex) polyhedra.\vspace*{-0.05in}

\begin{proposition}\label{prop.theo1.polyhedralGamma}
In the setting of Theorem~{\rm\ref{theorem.singular.constraint.mapping}}, suppose that the mapping $\varGamma$ is polyhedral. Then conditions {\rm(ii), (iv)}, and {\rm(v)} of Theorem~{\rm\ref{theorem.subderivative.singular.I}} are satisfied.
\end{proposition}\vspace*{-0.05in}
{\bf Proof}.
It is well known since \cite{Rob1981} (see also \cite[Section~3D]{Rock_ImpF}) that $\varGamma$ being polyhedral ensures its (uniform) outer Lipschitz continuity, i.e., the existence of $c\geq 0$ such that for any $(\xi,z)\in\dom\,\varGamma$, we can find $\varepsilon>0$ with
\begin{align*}
\varGamma(\xi',z')\subset \varGamma (\xi,z)+ c\Vert(\xi',z')-(\xi,z)\Vert\B\;\mbox{ whenever }
(\xi',z')\in (\xi,z)+\varepsilon\B.
\end{align*}
This yields condition (v) of Theorem~\ref{theorem.subderivative.singular.I}, while condition (ii) follows, e.g., from \cite[Proposition~5.12~(a)]{Rocka}. To verify (iv), we  check first that the mapping $H:=\varGamma(\bar{\xi},\cdot)$ is polyhedral. To proceed, observe that the polyhedrality of $\varGamma$ requires that
\begin{align*}
\gph\,\varGamma=\bigcup_{i=1}^s\left\{ (\xi,z,\eta)\mid A^i\xi+B^iz+C^i\eta\leq q^i\right\}
\end{align*}
for some matrices $A^1,\ldots,A^s$, $B^1,\ldots,B^s$, $C^1,\ldots, C^s$ and vectors $q^1,\ldots, q^s$. It is clear that we have $(z,\eta)\in\gph\,H=\gph\,\varGamma(\bar{\xi},\cdot)$ if and only if $\eta\in H(z)=\varGamma(\bar{\xi},z)$, i.e., $(\bar{\xi},z,\eta)\in\gph\,\varGamma$. Using the representation of $\gph\,\varGamma$ above gives us
\begin{align*}
\gph\,H=\bigcup_{i=1}^s\left\{ (z,\eta)\mid B^iz+C^i\eta\leq q^i-A^i\bar{\xi}\right\},
\end{align*}
which confirms that $H$ is a polyhedral mapping. Furthermore, the latter yields
\begin{align*}
\dom\,H & %=\{ z\mid H(z)\neq\varnothing \} & 
=\left\{ z\left|\,\exists i\in\{ 1,\ldots,s\}\,\exists\eta\;\mbox{ with }\;\; B^iz\leq q^i-(A^i\bar{\xi}+C^i\eta)\right.\right\}\\
~ & =\bigcup_{i=1}^s\left\{ z\left|\;(B^iz)_j\leq (q^i-A^i\bar{\xi})_j\;\mbox{ for all }\;j\;\mbox{ with }\;(C^i)_j=0^\top\right.\right\},
\end{align*}
where we use the fact that $\eta$ can be chosen arbitrarily large in norm when $(C^i)_j\neq 0$, and thus the upper bound $(q^i-(A^i\bar{\xi}+C^i\eta))_j$ of $(B^iz)_j$ can become arbitrarily large. This justifies the polyhedrality of $\dom\,H$.

To complete the verification of (iv), we represent $\dom\,H$ as the finite union of convex polyhedral sets $P^1,\ldots, P^s$ in order to get
\begin{align*}
\mathcal{A}(\bar{\xi})\dom H=\{\mathcal{A}(\bar{\xi}) z\mid \exists i\in\{ 1,\ldots,s\}\;\mbox{ with }\;z\in P^i\}=\bigcup_{i=1}^s\mathcal{A}(\bar{\xi})P^i.
\end{align*}
It follows from \cite[Proposition~3.55~(a)]{Rocka} that the set $\mathcal{A}(\bar{\xi})P^i$ is convex polyhedral for any $i\in\{ 1,\ldots,s\}$. Thus it is closed, which confirms (iv).\qed\vspace*{0.05in}

The next result can be extracted from \cite[Theorem~4.5]{JDiss24}. Its proof is similar and even simpler than that of Theorem~\ref{theorem.subderivative.singular.I}.\vspace*{-0.05in}

\begin{theorem}\label{theorem.subderivative.singular.II}
In the setting of Theorem~{\rm\ref{theorem.singular.constraint.mapping}}, let $\bar{\xi}\in \D$ satisfy conditions {\rm(i)--(iii)} of Theorem~{\rm\ref{theorem.subderivative.singular.I}} for some $\omega\in\R^n\times\R^p$. Impose also the following assumptions:
\begin{enumerate}
\item[\bf(iv')] The mapping $\K:\D\times\R^m\rightrightarrows\R^q$, defined as $\K(\xi,z):=\cone\left(\Gamma(\xi,z)\right)$, is outer semicontinuous relative to $\D\times\bd\B$ at $(\bar{\xi},z)$ for any $z\in \Z$.
\item[\bf(v')] It holds that
\begin{align}\label{f.nonseparate.OSC}
\textnormal{rge}\,\A(\bar{\xi}) \bigcap \varTheta \subset\{ 0\}
\end{align}
with the set $\varTheta$ taken from \eqref{f.subderivative.Theta}.
\end{enumerate}
Then we have the estimate
\begin{align*}
d\ell_\varXi(\bar{\xi})(\omega)\geq\min_{z\in \ZN}\left( \di\left[ 0,\A'(\bar{\xi};\omega) z+\varTheta+\textnormal{rge}\,\A(\bar{\xi})\right]\right).
\end{align*}
\end{theorem}

Theorem~\ref{theorem.subderivative.singular.II} contains the list of conditions that are different from those in Theorem~\ref{theorem.subderivative.singular.I}. Furthermore, the lower estimate for the subderivative of the LSV function in the latter theorem is also different from the former one. As follows from the discussions above, the combination of conditions (i) and (v') implies that Theorem~\ref{theorem.subderivative.singular.II} is applicable only in  rather special $\varXi$-singular points. 
Moreover, verifying whether condition (iv') of the theorem holds can be a challenging task. For this reason, we introduce below new conditions sufficient for (iv').\vspace*{-0.05in}

\begin{proposition}\label{prop1.suff.LSV.subderivative.II}
In the setting of Theorem~{\rm\ref{theorem.singular.constraint.mapping}}, suppose that $\varGamma$ is outer semicontinuous and cone-valued. Then condition {\rm(iv')} of Theorem~{\rm\ref{theorem.subderivative.singular.II}} is fulfilled.
\end{proposition}\vspace*{-0.05in}
{\bf Proof}. We have $\mathcal{K}(\xi,z)=\cone(\varGamma(\xi,z))=\varGamma(\xi,z)$, which tells us that the outer semicontinuity of $\varGamma$ yields the one for $\mathcal{K}$, and thus condition (iv') is valid.\qed\vspace*{-0.05in}

\begin{proposition}\label{prop2.suff.LSV.subderivative.II}
In the setting of Theorem~{\rm\ref{theorem.singular.constraint.mapping}}, suppose that $\varGamma$ is outer semicontinuous and that the condition
\begin{align*}
H(z)=\{ 0\}\times\mathcal{T}(z)\;\mbox{ whenever }\;z\in\R^m
\end{align*}
holds for $\mathcal{T}:\R^m\rightrightarrows\R^{q-s}$ with $s\leq q$. Then the following two requirements held together ensure the fulfillment of condition {\rm(iv')} of Theorem~\ref{theorem.subderivative.singular.II}:
\begin{align*}
\mathcal{T}(0)=\{ 0\},\qquad \mathcal{T}^{-1}(0)=\{ 0\}.
\end{align*}
\end{proposition}
{\bf Proof}. Consider the sequences $\{ (\xi^k,z^k)\}\subset \mathcal{D}\times\bd\B$ and $\{ \eta^k\}$ with $\eta^k\in\mathcal{K}(\xi^k,z^k)$ converging to $(\bar{\xi},z)$ for some $z\in\Z=\bd\B \cap\dom\,H$ and $\eta$, respectively. We aim to show that $\eta\in\mathcal{K}(\bar{\xi},z)$. If $\eta=0$, then the inclusion $\eta=0\in \mathcal{K}(\bar{\xi},z)$ is trivial. Assuming now that $\eta\neq 0$ tells us by the definition of $\mathcal{K}$ that for any $k$ there are $\gamma_k\geq 0$ and $\zeta^k\in\varGamma(\xi^k,z^k)$ satisfying $\eta^k=\gamma_k\zeta^k$. If $\gamma_k\rightarrow\infty$, then we can assume that $\gamma_k\neq 0$ for all $k$ and thus get $\gamma_k^{-1}\eta^k=\zeta^k\rightarrow 0$. It follows from the outer semicontinuity of $\varGamma$ and the special structure of $H$ that $0\in\varGamma(\bar{\xi},z)=H(z)=\{ 0\}\times\mathcal{T}(z)$, which implies by the imposed assumption that $\mathcal{T}^{-1}(0)=\{ 0\}$, a contradiction.
	
It remains to examine the case where $\{\gamma_k\}$ is bounded. If in this case we have $\gamma_k\rightarrow 0$ along a subsequence, then the condition $\eta\neq 0$ yields $\Vert\zeta^k\Vert\rightarrow\infty$. Since $\varGamma$ is positively homogeneous in the second argument, we find $\eta^k=\gamma_k\zeta^k\in \varGamma(\xi^k,\gamma_kz^k)$ for all $k$, hence, $\eta\in \varGamma(\bar{\xi},0)=H(0)=\{ 0\}\times\mathcal{T}(0)$ by the outer semicontinuity of $\varGamma$, $\gamma_k\rightarrow0$, and the special structure of $H$. The imposed assumptions ensure that $H(0)=\{ 0\}\times\mathcal{T}(0)=\{ 0\}$  and hence $\eta=0$, a contradiction. Therefore, both sequences $\{ \gamma_k\}$ and $\{\zeta^k\}$ are bounded, which yields their convergence along some subsequences. Denoting by $\gamma\ge 0$ and $\zeta\in\R^q$ their limiting points gives us $\eta=\gamma\zeta$. It follows from the outer semicontinuity of $\varGamma$ that $\zeta\in\varGamma(\bar{\xi},z)$, which brings us to $\eta=\gamma\zeta\in \mathcal{K}(\bar{\xi},z)$ and thus confirms that $\mathcal{K}$ is outer semicontinuous.\qed\vspace*{0.05in}

We conclude this section by three examples illustrating relationships between the assumptions and conclusions of the obtained results. In particular, Example~\ref{example.2} shows that the assumptions of Theorem~\ref{theorem.subderivative.singular.I} may be fulfilled, while those in Theorem~\ref{theorem.subderivative.singular.II} are not. Example~\ref{example.3} illustrates the opposite.\vspace*{-0.05in}

\begin{example}\label{example.2}
Let $\D:=\R$, $\A(\cdot)\equiv(1,0)$, $\bar{\xi}:=0$, and
\begin{align*}
\Gamma(\xi,z_1,z_2):=\left\{ 
\begin{array}{c l}
\varnothing & \textnormal{ if }\xi\neq 0 \textnormal{ or }(z_1,z_2)\notin\{ 0\}\times\R_{-},\\
\{ 0,-1\} & \textnormal{ if }\xi= 0 \textnormal{ and }(z_1,z_2)\in\{ 0\}\times\R_{-}.
\end{array}\right.
\end{align*}
The mapping $\varGamma$ is clearly polyhedral, so Proposition~\ref{prop.theo1.polyhedralGamma} guarantees the fulfillment of (ii), (iv), and (v) in Theorem~\ref{theorem.subderivative.singular.I}.
In particular, it can be seen that condition (v) therein holds with $c=0$. Moreover, conditions (i) and (iii) in the theorem hold as well, and we have $\dom\,H=\{ 0\}\times\R_{-}$ and $\A(\bar{\xi})\dom\,H=\{ 0\}$. Therefore, Theorem~\ref{theorem.subderivative.singular.I} is applicable. 
At the same time, we get $\textnormal{rge}\,\mathcal{A}(\bar{\xi})=\R$ and $\varTheta= \R_{-}$, which yield $\textnormal{rge}\,\A(\bar{\xi})\cap\varTheta=\R_{-}$, i.e., condition (v') of Theorem~\ref{theorem.subderivative.singular.II} fails.\qed
\end{example}\vspace*{-0.15in}

\begin{example}\label{example.3}
For a nonempty closed set $\mathcal{D}\subset\R^n$, let $\mathcal{A}\equiv 0$ and define $\varGamma(\xi,z):=\mathcal{N}(\xi)$ with a nonempty cone-valued outer semicontinuous mapping $\mathcal{N}:\mathcal{D}\rightrightarrows\R^q$. Here we have $\K=\varGamma$ and can easily check that for any $\bar{\xi}\in\D$ all the assumptions of Theorem~\ref{theorem.subderivative.singular.II} are fulfilled. At the same time, condition (v) in Theorem~\ref{theorem.subderivative.singular.I} may not hold without further specifications of $\mathcal{N}$, i.e., the theorem is not applicable without imposing additional assumptions.\qed
\end{example}\vspace*{-0.05in}
	
The final example shows that condition (iv') of Theorem~\ref{theorem.subderivative.singular.II} cannot be removed from the list of assumptions without destroying the conclusion of the theorem.\vspace*{-0.05in}

\begin{example}\label{example.4}
Let $\D:=\R$, $\bar{\xi}:=0$, $\omega:=1$, $\A(\xi):=(\xi,\xi^2)$, and
\begin{align*}
\varGamma(\xi,z):=\left\{ 
\begin{array}{c l}
\varnothing & \textnormal{ if }z\notin\R\times\{ 0\},\\
\{-\xi\} & \textnormal{ if }z\in\R\times\{ 0\}.
\end{array}
\right.
\end{align*}
It is straightforward to confirm the fulfillment of  (i)--(iii), and (v') of  Theorem~\ref{theorem.subderivative.singular.II}. Furthermore, it follows that $\ell_{\varXi}(\bar{\xi}+t\omega)=0$ for all $t\geq 0$, and so the equality $d\ell_{\varXi}(\bar{\xi})(\omega)=0$ holds. However, we have $\mathcal{Z}_0=\{(-1,0),(1,0)\}$, $\rge\mathcal{A}(\bar{\xi})=\{ 0\}$, $\mathcal{A}'(\bar{\xi};\omega)=(1,0)$, and $\varTheta=\{ 0\}$. Therefore, the conclusion of Theorem~\ref{theorem.subderivative.singular.II} fails, and the only possible reason for this is the violation of condition (iv').\qed
\end{example}\vspace*{-0.32in}

\section{Metric Regularity via Coderivatives and the LSV Function}\label{sec.MR}\vspace*{-0.1in}

The preceding section dealt with with the study of the LSV function by using the subderivative, which is a {\em primal-space} construction of generalized differentiation. On the other hand, it has been well recognized that {\em dual-space} constructions revolving around {\em co}derivatives play a prominent role in characterizing metric regularity and related well-posedness properties of variational analysis. To incorporate  such constructions in our study of mixed-order regularity, this section reviews appropriate dual-space notions of generalized differentiation and unify them with the LSV function in understanding the classical well-posedness properties for structural set-valued mappings. We begin with the fundamental construction of the (limiting) {\em normal cone} to closed sets introduced (in the equivalent form) by the second author \cite{Mo76}.\vspace*{-0.1in}  

\begin{definition}\label{def.NC} For any set $C_0\subset\R^n$ with $\bar x\in C_0$, the $($limiting, basic, Mordukhovich$)$ {\em normal cone} to $C_0$ at $\bar x$ is given by
\begin{align}\label{nc}
N_{C_0}(\bar x):=\Limsup_{x\overset{C_0}{\rightarrow}\bar{x}}\left( T_{C_0}(x)^\circ\right),
\end{align}
where $T_{C_0}(x)^\circ:=\{ v\mid v^\top w\leq 0\;\mbox{ for all }\;w\in T_{C_0}(x)\}$ is the polar cone of $T_{C_0}(x)$. We define $N_{C_0}(\bar{x}):=\emptyset$ when $\bar x\notin C_0$.
\end{definition}\vspace*{-0.07in}

Various equivalent representations and properties of the normal cone \eqref{nc} can be found in the books \cite{Mo18,Mo24,Rocka} and the references therein. For convex sets $C_0$, the normal cone \eqref{nc} agrees with the classical normal cone of convex analysis, while the values of \eqref{nc} may not be convex even for simple nonconvex sets like $C_0={\rm epi}(-|x|)$ and $C_0=\gph(|x|)$ at $(0,0)$. Nevertheless, this normal cone and the associated generalized differential construction for extended-real-valued functions and set-valued mappings (subdifferential and coderivative) enjoy {\em full calculi} based on {\em variational/extremal principles} of variational analysis; see the aforementioned books for more details.\vspace*{0.05in}

In this paper, we use two simple properties of the normal cone \eqref{nc}, which follow directly directly from Definition~\ref{def.NC}.\vspace*{-0.1in}

\begin{lemma}\label{l.NC.osc}
For any closed set $C_0\subset\R^n$, the normal cone mapping $N_{C_0}:\R^n\rightrightarrows\R^n$ associated with $C_0$ is outer semicontinuous.
\end{lemma}\vspace*{-0.2in}

\begin{lemma}\label{l.NC.product}
For any sets $C_1\subset\R^{p}$ and $C_2\subset\R^q$ with $(x,z)\in C_0:=C_1\times C_2$, we have
\begin{align*}
N_{C_0}(x,z)=N_{C_1}(x)\times N_{C_2}(z).
\end{align*}
\end{lemma}\vspace*{-0.05in}

The main dual-space construction used for characterizing metric regularity and related well-posedness properties is the following {\em coderivative} of set-valued mapping introduced in \cite{Mo80} by using the normal cone \eqref{nc} to the mapping graph.\vspace*{-0.07in}

\begin{definition}\label{def.coderivatives}
Given a set-valued mapping $S:\R^n\rightrightarrows\R^m$ and a point $(\bar x,\bar y)\in\gph\, S$,
the \textit{coderivative} of $S$ at $(\bar x,\bar y)$ is the set-valued mapping $D^*S(\bar x|\bar y):\R^m\rightrightarrows\R^n$ defined by
\begin{align}\label{cod}
D^*S(\bar x|\bar y)(z):=\left\{ v\in\R^n\;\left|\; (v,-z)\in N_{\gph\,S}(\bar x,\bar y)\right.\right\}.
\end{align}
\end{definition}

Based on the definition and Lemma~\ref{l.NC.osc}, we have the following properties.\vspace*{-0.1in}

\begin{lemma}\label{l.coder.osc}
Given a mapping $S:\R^n\rightrightarrows\R^m$ and a point $(\bar x,\bar y)\in\gph\,S$ around which $\gph\,S$ is closed, the coderivative \eqref{cod} is positively homogeneous and outer semicontinuous with respect to all its variables.
\end{lemma}\vspace*{-0.05in}

For a set-valued mapping $S:\R^n\rightrightarrows\R^m$ with the closed graph $\D:=\gph\,S$, define the new mapping $\varXi:\gph\,S\times\R^m\rightrightarrows\R^n$ by
\begin{align*}
((u,y),z)\mapsto\varXi((u,y),z):=D^*S(u|y)(z)
\end{align*}
and observe by Lemma~\ref{l.coder.osc} that $\varXi$ satisfies the requirements in Definition~\ref{def.LSV}. Thus the LSV function \eqref{LSV} for $\varXi$ is well-defined, and from now on we can use the notation:
\begin{align}\label{f.coderi.LSV}
\Reg(u,y;S):=\ell_{\varXi}(u,y).
\end{align}
Lemma~\ref{l.coder.osc} tells us that the mapping $\varXi$ is outer semicontinuous, and so it follows from Theorem~\ref{theorem.singular.constraint.mapping} that a point $(u,y)\in\gph\,S$ is {\em $\varXi$-singular if and only if}
\begin{align}\label{f.Cod.singular}
\exists\,z\neq 0\;\mbox{ with }\;0\in D^*S(u|y)(z),
\end{align}
i.e., the coderivative of $S$ at $(u,y)$ is {\em singular}.\vspace*{0.02in}

To proceed further, we need the following three lemmas about the coderivative properties. The first lemma is the coderivative sum rule taken from \cite[Theorem~1.62(ii)]{Mo06}, which holds for mappings between general normed spaces and where the graph closedness is not required.\vspace*{-0.07in}

\begin{lemma}\label{lemma.coderivative.sum} Let $F:\R^n\rightarrow\R^m$ be strictly differentiable at $\bar x\in\R^n$, and let $\bar y\in\R^m$ be such that $(\bar x,\bar y)\in\gph\,C$ for some set-valued mapping $C:\R^n\rightrightarrows\R^m$. Then
\begin{align*}
D^*(F+C)(\bar x|F(\bar x)+\bar y)(z)=\nabla F(\bar x) z+D^*C(\bar  x|\bar y)(z)\;\mbox{ for all }\;z\in\R^m.
\end{align*}
In particular, we have the LSV function values calculated by
\begin{align*}
\Reg(\bar x,F(\bar x)+\bar y;F+C)=\inf_{\Vert z\Vert=1} \di\left[ 0,\nabla F(\bar x)z+D^*C(\bar x|\bar y)(z)\right].
\end{align*}
\end{lemma}\vspace*{-0.05in}
 
The second lemma allows us to deal with the coderivative of product mappings. It is a direct consequence of the coderivative definition \eqref{cod} and Lemma~\ref{l.NC.product}. \vspace*{-0.05in}

\begin{lemma}\label{lemma.coderivative.product} Given $R:\R^k\rightrightarrows\R^p$ and $T:\R^l\rightrightarrows\R^q$, consider the product mapping $\mathcal{P}:\R^k\times\R^l\rightrightarrows\R^{p+q}$ defined by
\begin{align}\label{f.product.map}
\mathcal{P}(x,u):=\left\{\left.\left(\begin{array}{c}
y\\v\end{array}			
\right)\in\R^{p+q}\right|\, y\in R(x),\;v\in T(u)\right\}.
\end{align}
Then for all $\left((\bar x,\bar u),(\bar y,\bar v)\right)\in\gph\mathcal{P}$, and $(\theta,\zeta)\in\R^{p+q}$, we have
\begin{align*}
D^*\mathcal{P}\left((\bar x,\bar u)|(\bar y,\bar v)\right)\left(\begin{array}{c}
\theta\\\zeta\end{array}			
\right) & =\left\{\left.(v,\mu)\in\R^k\times\R^l\right|\, \begin{array}{r l}
v & \in D^*R(\bar x|\bar u)(\theta),\\ 
\mu & \in D^*T(\bar y|\bar v)(\zeta)
\end{array}\right\}.
\end{align*}
\end{lemma}\vspace*{0.05in}

The third lemma concerns the coderivative of the {\em indicator mapping} associated with a given set $\Omega$ by the correspondence
\begin{align}\label{indic}
\Delta_{\Omega}(x):=\left\{ 
\begin{array}{c l}
\{ 0\} & \textnormal{ if }\;x\in\Omega,\\
\varnothing & \textnormal{ if }\;x\notin\Omega.
\end{array}
\right.
\end{align}
The proof follows from \eqref{cod} due to $\gph\,\Delta_\Omega=\Omega\times\{0\}$; see \cite[Proposition~1.33]{Mo06}.\vspace*{-0.05in}

\begin{lemma}\label{prop.coderivative.calculus}
Given $C_0\subset \R^q$, $C_1\subset \R^p$ and $\bar x\in C_0$, $\bar y\in C_1$, we have
\begin{align*}
D^*\varDelta_{C_0}(\bar x|0)(z)=N_{C_0}(\bar x)\;\mbox{ for all }\;z\in \R^q.
\end{align*}
Furthermore, it holds for the constant mapping $M(x)\equiv C_1$ that
\begin{align*}
D^*M(\bar x|\bar y)(u)=\varDelta_{N_{C_1}(\bar y)}(-u)\;\mbox{ whenever }\;u\in\R^p.
\end{align*}
\end{lemma}\vspace*{-0.05in}

Now we recall the fundamental notions of {\em well-posedness} in variational analysis. The terminology ``around the point" is used to highlight that the corresponding properties provide {\em robust} behavior in a {\em neighborhood} of the point in question.\vspace*{-0.05in}

\begin{definition}\label{definition.metric.regularity}
A set-valued mapping $S:\R^n\rightrightarrows\R^m$ is \textit{metrically regular} around $(\bar u,\bar y)\in\gph\,S$ with a constant $c> 0$ if there exists $\varepsilon>0$ such that
\begin{align*}
\di\left[ u,S^{-1}(y)\right]\leq c\cdot \di\left[ y, S(u)\right]\;\mbox{ for all }\;(u,y)\in (\bar u,\bar y)+\varepsilon\B.
\end{align*}
The mapping $S$ is \textit{Lipschitz-like} around $(\bar u,\bar y)$ with a constant $c> 0$ if there exists
$\varepsilon>0$ for which we have the inclusion
\begin{align*}
S(u')\cap(\bar y+\varepsilon\B)\subset S(u)+c\Vert u'-u\Vert\B\;\mbox{ whenever }\;u,u
'\in u+\varepsilon\B.
\end{align*}
\end{definition}\vspace*{0.05in}
The Lipschitz-like property is also known as the {\em Aubin $($pseudo-Lipschitz$)$} property. It is actually a graphical localization of the classical (Hausdorff) local Lipschitz continuity for set-valued mappings around points of their domains, and so we choose to emphasize the Lipschitzian nature in the name of the above graphical localization; see the books \cite{Mo06,Rocka} for more references and discussions.
\vspace*{0.03in}

The following coderivative characterizations of the well-posedness properties from Definition~\ref{definition.metric.regularity} goes back to \cite{Mo93} being called the {\em Mordukhovich criterion} in \cite{Rocka}. It plays a fundamental role in variational analysis and its applications; see, e.g., the books \cite{Mo06,Mo18,Mo24,thibault} with the vast bibliographies  therein.\vspace*{-0.05in}

\begin{theorem}\label{theorem.metric-regularity}
Let $S:\R^n\rightrightarrows\R^m$ be a set-valued mapping whose graph is closed around $(\bar u,\bar y)\in\gph\,S$.
Then we have the equivalent statements:\vspace*{-0.05in}
\begin{enumerate}
\item[\bf(a)] $S$ is metrically regular around $(\bar u.\bar y)$ with some constant $c>0$.

\item[\bf(b)] $S^{-1}$ is Lipschitz-like around $(\bar y,\bar u)$ with some constant $\ell>0$. 

\item[\bf(c)] There are constants $\theta,\kappa>0$ such that 
\begin{align}\label{f.linear.open}
S(u)\cap (\bar y+\theta\B)+\tau\B\subset S(u+\kappa\tau\B)
\end{align}
for all $(u,\tau)\in\R^n\times\R_+$ with $\Vert u-\bar u\Vert+\tau\leq\theta$.
\item[\bf(d)] The coderivative of $S$ at $(\bar u,\bar y)$ is not singular meaning that \eqref{f.Cod.singular} does not 
hold, and thus ${\rm ker}\,D^*S(\bar u,\bar y)=\{0\}$.
\end{enumerate}\vspace*{-0.07in}
Moreover, the infimum of all $c>0$ satisfying {\rm(a)} agrees with the exact bounds of constants $\ell$ and $\kappa$ in {\rm(b)} and {\rm(c)} being equal to $\Reg(\bar u,\bar y;S)^{-1}$.
\end{theorem}\vspace*{-0.05in}

Our next result provides relations between LSV function values for closed-graph mappings represented as sums of smooth single-valued and set-valued terms.\vspace*{-0.05in} 

\begin{theorem}\label{theorem.metric.regularity.QCM}
Let $F:\R^n\rightarrow\R^m$ be continuous, and let $C:\R^n\rightrightarrows\R^m$ be closed-graph. Define the mappings $\varSigma:\R^n\rightrightarrows\R^m$, $\mathcal{C}:\R^n\rightrightarrows\R^n\times\R^m$, $Q:\R^n\times\R^m\rightrightarrows\R^m$, and $\mathcal{Q}:\R^n\times\R^n\times\R^m\rightrightarrows\R^m$ by, respectively,\vspace*{-0.1in}
\begin{align*}
\varSigma(u) & :=F(u)+C(u),\\
\mathcal{C}(u)& :=\left(
\begin{array}{c}
-u\\
F(u)
\end{array}\right)+\gph\;C,\\
Q(u,y) & :=(F(u)+y)+\varDelta_{\gph\,C}(u,y),\\
\mathcal{Q}(u,\sigma,y) & :=\left( \begin{array}{c}
-u+\sigma\\
F(u)+y
\end{array}
\right)+\varDelta_{\gph\,C}(\sigma,y).
\end{align*}\vspace*{-0.1in}
Pick $(\bar u^,\bar y)\in\gph\,C$ and assume that
$F$ is strictly differentiable at $\bar u$. Then the numbers\vspace*{0.03in}
\begin{align*}
R_\varSigma & :=\Reg(\bar u,F(\bar u)+\bar y;\varSigma),\\ 
R_{\mathcal{C}} & :=\Reg(\bar u,(0,F(\bar u)+\bar y);\mathcal{C}),\\
R_Q & :=\Reg((\bar u,\bar y),F(\bar u)+\bar y;Q),\\
R_{\mathcal{Q}} & :=\Reg((\bar u,\bar u,\bar y),(0,F(\bar u)+\bar y);\mathcal{Q})
\end{align*}
are well-defined, and satisfy the relationships
$R_\varSigma \geq \max\{ R_\mathcal{C},R_Q\}\geq R_\mathcal{Q}
\geq 0$. Moreover, the equality $R_\mathcal{Q}=0$ implies that $R_\varSigma=R_\mathcal{C}=R_Q=0$.
\end{theorem}\vspace*{-0.02in}
{\bf Proof}. The claim that the aforementioned numbers are well-defined is a consequence of the fact that the corresponding mappings are closed-graph, which follows from the assumptions of the theorem. Indeed, this fact allows us to represent values of the LSV function \eqref{LSV} in the coderivative form \eqref{f.coderi.LSV}. Using now Lemma~\ref{lemma.coderivative.sum} and Lemma~\ref{prop.coderivative.calculus} yields the coderivative representations
\begin{align*}
D^*\varSigma(\bar u^*|F(\bar u)+\bar y)(z ) & =\nabla F(\bar u) z+D^*C(\bar u|\bar y)(z),\\
D^*\mathcal{C}(\bar u|(0,F(\bar u)+\bar y))(v,z) & =-v+\nabla F(\bar u) z+\varDelta_{N_{\gph\,C}(\bar u,\bar y)}(-(v,z)),\\
D^*Q((\bar u,\bar y)|F(\bar u)+\bar y)(v,z) & =\left(\begin{array}{c}
\nabla F(\bar u) z\\
z\end{array}\right)
+N_{\gph\,C}(\bar u,\bar y),\\
D^*\mathcal{Q}((\bar u,\bar u,\bar y)|(0,F(\bar u)+\bar y))(v,z) & =\left(\begin{array}{c}
-v+\nabla F(\bar u) z\\
v\\z\end{array}	\right)+\{ 0\}\times N_{\gph\,C}(\bar u,\bar y).
\end{align*}
To verify the estimate $R_\varSigma \geq \max\{ R_\mathcal{C},R_Q\}$, we get from the latter formulas and the definitions of $R_\mathcal{C},\; R_\varSigma$ and of the coderivative \eqref{cod} that
\begin{align*}
R_\mathcal{C} & =\inf_{\Vert (v,z)\Vert=1}\di\left[0,D^*\mathcal{C}(\bar u|(0,F(\bar u)+\bar y))(v,z)\right]\\
~ & =\inf_{\substack{\Vert (v,z)\Vert=1\\ v\in D^*C(\bar u|\bar y)(z)}}\left\Vert v+\nabla F(\bar u) z\right\Vert\\
~ & \leq\inf_{\Vert z\Vert=1}\di\left[0,\nabla F(\bar u) z+D^*C(\bar u|\bar y)(z)\right]=R_\varSigma.	
\end{align*}
In a similar way, we get  $R_Q\leq R_\varSigma$, and hence $\max\{ R_\mathcal{C},R_Q\}\leq R_\varSigma$. The second estimate $\max\{ R_\mathcal{C},R_Q\}\geq R_\mathcal{Q}$ follows from the relationships
\begin{align*}
R_\mathcal{Q} & =\inf_{\Vert (v,z)\Vert=1}\di\left[0,D^*\mathcal{Q}((\bar u,\bar u,\bar y)|(0,F(\bar u)+\bar y))(v,z)\right]\\
~ & =\inf_{\Vert (v,z)\Vert=1}\left\{ \left.\left\Vert \left(
\begin{array}{c}
-v+\nabla F(\bar u) z\\
v+\chi\\
z+\zeta
\end{array}
\right)\right\Vert\right|\;(\chi,\zeta)\in N_{\gph\,C}(\bar u,\bar y)   \right\}\\~ & \leq\inf_{\substack{\Vert (v,z)\Vert=1\\(v,z)\in N_{\gph\,C}(\bar u,\bar y)}}\left\Vert -v+\nabla F(\bar u) z\right\Vert=R_\mathcal{C}\leq \max\{ R_{\mathcal{C}},R_Q\}.
\end{align*}
It remains to justify the last claim of the theorem. Assume that $R_\mathcal{Q}=0$. By the coderivative criterion of Theorem~\ref{theorem.metric-regularity}, this means that
\begin{align*}
\exists\,(\bar{v},\bar{z})\neq(0,0)\;\mbox{ such that }\;(\bar{v},-\bar{z})\in N_{\gph\,C}(\bar u,\bar y),\quad \nabla F(\bar u) \bar{z}+\bar{v}=0.
\end{align*}
For the point $(\bar{v},\bar{z})$ satisfying the latter, we see that $\bar{z}\neq 0$ (otherwise, it follows from $\bar{z}=0$ that $\bar{v}=0$, a contradiction). This gives us $0\in \nabla F(\bar u) \bar{z}+D^*C(\bar u|\bar y)(\bar{z})$ with $\bar{z}\neq 0$, which yields $0=R_\varSigma\geq \max\{ R_\mathcal{C},R_Q\}\geq 0$ and thus completes the proof. \qed
\vspace*{-0.2in}

\section{Criteria for Metric 
2-Regularity}\label{sec.SC}\vspace*{-0.1in}

This section is devoted to the notion of {\em metric} 2-{\em regularity} for set-valued mappings that constitutes another well-posedness concept of variational analysis, which works when the standard  metric regularity fails. Observe that the established metric regularity criteria concern the case where $\textnormal{Reg}(\bar u,\bar y;S)>0$. From now on, we investigate the opposite case, i.e., when $\textnormal{Reg}(\bar u,\bar y;S)=0$, which is precisely the setting where $S$ is {\em not metrically regular} around $(\bar u,\bar y)$.

For a point $\bar u\in\R^n$ and a direction $w\in\R^n$, consider the {\em directional neighborhood}
\begin{align}\label{f.dir.neigh}
K_{\varepsilon,\delta}(\bar u;w):=\bar u+\left(\varepsilon\B\cap \cone(w+\delta\B)\right),
\end{align}
where $\varepsilon,\delta>0$. It is clear that $K_{\varepsilon,\delta}(\bar u;0)=\bar u+\varepsilon\B$ and  $K_{\varepsilon,\delta}(\bar u;r\cdot w)=K_{\varepsilon, r^{-1}\delta}(\bar u;w)$ for all $r>0$.
Along with the coderivative, we explore here yet another common generalized first-order derivative of the primal type; see \cite{Rocka} for more details.\vspace*{-0.05in}

\begin{definition}\label{def.graph.derivative}
The \textit{graphical derivative} of $S:\R^n\rightrightarrows\R^m$ at $(\bar u,\bar y)\in\gph\,S$ is the set-valued mapping $DS(\bar u|\bar y):\R^n\rightrightarrows\R^m$ defined by
\begin{align*}
DS(\bar u|\bar y)(w):=\left\{ v\in\R^n\;\left|\; (w,v)\in T_{\gph\,S}(\bar u,\bar y)\right.\right\},\quad w\in\R^m,
\end{align*}
via the tangent cone \eqref{tan} to the mapping graph.
\end{definition}\vspace*{-0.05in}

Now we recall another concept of regularity, recently coined in \cite{JDiss24}, which we will work with in this section.\vspace*{-0.05in}

\begin{definition}\label{definition.gamma.pre.regularity}
A mapping $S:\R^n\rightrightarrows\R^m$ is \textit{metrically $2$-regular} around $(\bar u,\bar y)\in\gph\,S$ relative to a direction $w\neq 0$ with constant $c>0$ if there are $\varepsilon_0,\delta_0,\rho_0>0$ with $\rho_0>\Vert w\Vert^{-1}\sup\{ \Vert\eta\Vert\mid \eta\in DS(\bar u|\bar y)(w)\}$ such that $S$ is metrically regular around all $u\in K_{\varepsilon_0,\delta_0}(\bar u;w)\setminus\{\bar u\}$ for any $y\in S(u)\cap (\bar y+\rho_0\Vert u-\bar u\Vert\B)$ with constant $c/\Vert u-\bar u\Vert$.
\end{definition}\vspace*{-0.05in}

Since the graphical derivative values are closed sets, we see that the metric 2-regularity of 
$S$ around $(\bar u,\bar y)$ relative to $w\neq 0$ ensures that $DS(\bar u|\bar y)(w)\neq\varnothing$.\vspace*{-0.05in}

\begin{remark}\label{rem.m2r.mr}
If $S$ is metrically regular around $(\bar u,\bar y)$, then it is metrically 2-regular around this point in any nonzero direction $w\in\dom(DS(\bar u|\bar y))$. The converse is not necessarily true as shown by the example where $S(u):=u^2$ for $u\in\R$ and $\bar u:=0$.\qed
\end{remark}\vspace*{-0.2in}

\begin{remark}\label{rem.directional}
Metric 2-regularity is a {\em directional} notion, which is not defined for the zero vector. Furthermore, if $S$ is metrically 2-regular around $(\bar u,\bar y)$ relative to $w$, then there exist sequences $t_k\searrow 0$, $w^k\rightarrow w$, and $y^k\in S(\bar u+t_k\bar w^k)$ with $\Vert y^k-\bar y\Vert\leq \rho_0t_k\Vert w^k\Vert$ such that $S$ is metrically regular around $(\bar u+t_kw^k,y^k)$ with the constant $c/(t_k\Vert w^k\Vert)$ for all $k$. Thus this property yields the existence of a sequence of points from the mapping graph around which it is metrically regular at a special rate, possibly tending to infinity. This is in contrast to metric regularity, which is maintained around all neighboring points from the graph with the same regularity-mode. \qed
\end{remark}\vspace*{-0.1in}

Observe that the metric 2-regularity relative to $w\neq 0$ holds if and only if it holds relative to $w/\Vert w\Vert$. Therefore, it is sufficient to consider only unit directions in the context of metric 2-regularity.\vspace*{0.05in}

The next lemma actually rephrases Theorem~\ref{theorem.metric-regularity} in the setting convenient for our study of metric 2-regularity; see a detailed proof in \cite{JDiss24}.\vspace*{-0.1in}

\begin{lemma}\label{lemma.gamma-regularity} Consider set-valued mappings $\varPsi,S:\R^n\rightrightarrows\R^m$ with $\gph\,\varPsi\subset\gph\,S$, a unit direction $w\in\R^n$, and a point $(\bar u,\bar y)\in\gph\,\varPsi$. If the set $\gph\,S$ is closed,  then the following statements are equivalent:\vspace*{-0.05in}
\begin{enumerate}
\item[\bf(a)] There exist positive numbers $\varepsilon,\delta,c$ such that $S$ is metrically regular around any $u\in K_{\varepsilon,\delta}(\bar u;w)\setminus\{\bar u,y)$ as $y\in\varPsi(u)$ with the constant $c/\Vert u-\bar u\Vert$.

\item[\bf(b)] There exist positive numbers $\varepsilon,\delta,c$ such that for all $y\in\varPsi\left(K_{\varepsilon,\delta}(\bar u;w)\right)$
and all $u\in \varPsi^{-1}(y)\cap K_{\varepsilon,\delta}(\bar u;w)\setminus\{\bar u\}$, the mapping $S^{-1}$ is Lipschitz-like around $(y,u)$ with the constant $c/\Vert u-\bar u\Vert$.

\item[\bf(c)] There exist positive numbers $\varepsilon,\delta,c$ such that for all $\widehat{u}\in K_{\varepsilon,\delta}(\bar u;w)$ and all $\widehat{y}\in\varPsi(\widehat{u})$, we find $\theta(\widehat{u},\widehat{y})>0$ for which the inclusion
\begin{align}\label{f.gamma-reg.covering}
S(u)\cap \left( \widehat{y}+\theta(\widehat{u},\widehat{y})\B\right)+\tau\Vert \widehat{u}-\bar u\Vert\B\subset S\left(u+c\tau\B\right)
\end{align}
is fulfilled whenever $(u,\tau)\in\R^n\times\R_+$ with $\Vert u-\widehat{u}\Vert+\tau\leq \theta(\widehat{u},\widehat{y})$.
		
\item[\bf(d)] There exist positive numbers $\varepsilon,\delta,c$ such that
\begin{align*}
c^{-1}\Vert u-\bar u\Vert\leq \inf_{y\in\varPsi(u)}\textnormal{Reg}(u,y;S)\;\mbox{ for all }\; u\in K_{\varepsilon,\delta}(\bar u;w).
\end{align*}
\end{enumerate}
The exact bounds of the corresponding constants $c>0$ in {\rm(a)--(d)} are the same.
\end{lemma}\vspace*{-0.05in}

With the above lemma at hand, we can formulate a counterpart of Theorem~\ref{theorem.metric-regularity} for the case of metric 2-regularity.\vspace*{-0.05in}

\begin{theorem}\label{theorem.m2r.criterion}
For a closed-graph mapping $S:\R^n\rightrightarrows\R^m$ and $(\bar u,\bar y)\in\gph\,S$, assume that $DS(\bar u|\bar y)(w)\neq\varnothing$ with some unit direction $w\in \R^n$. The following are equivalent:\vspace*{-0.05in}
\begin{itemize}
\item[\bf(a)] $S$ is metrically 2-regular around $(\bar u,\bar y)$ relative to $w$ with a constant $c>0$.

\item[\bf(b)] There is a number $\rho_0>\Vert w\Vert^{-1}\max\{ \Vert\eta\Vert\mid \eta\in DS(\bar u|\bar y)(w)\}$ such that assertions {\rm(a)--(d)} of Lemma~{\rm\ref{lemma.gamma-regularity}} are satisfied for the mapping $\varPsi:\R^n\rightrightarrows\R^m$ defined by
\begin{align*}
\Psi(u):=S(u)\cap(\bar y+\rho_0\Vert u-\bar u\Vert\B).
\end{align*}

\item[\bf(c)] We have the implication
\begin{align}\label{f.m2r.criterion}
\textnormal{Reg}(\bar u,\bar y;S)=0\Longrightarrow\inf_{\eta\in DS(\bar u|\bar y)(w)}d\textnormal{Reg}(\bar u,\bar y;S)(w,\eta)>0.
\end{align}
\end{itemize}
In particular, if $\textnormal{Reg}(\bar u,\bar y;S)=0$, then the infimum of all $c>0$ for which $S$ is metrically 2-regular around $(\bar u,\bar y)$ relative to $w$ is equal to
\begin{align*}
\sup\limits_{\eta\in DS(\bar u|\bar y)(w)}\left(\frac{1}{ d\textnormal{Reg}(\bar u,\bar y;S)(w,\eta)}\right).
\end{align*}
\end{theorem}
{\bf Proof}.
Equivalence (a)$\Longleftrightarrow$(b) follows from the definition of metric 2-regularity and Lemma~\ref{lemma.gamma-regularity}. The remaining claims are verified in \cite[Theorem~4.18]{JDiss24}.
\qed\vspace*{0.03in}

The next result provides a characterization of metric 2-regularity for mappings represented as sums of smooth single-valued and set-valued mappings. At this point, parallels to the structure of the mappings relevant to Theorems~\ref{theorem.singular.constraint.mapping}--\ref{theorem.subderivative.singular.II} become clear.\vspace*{-0.05in}

\begin{theorem}\label{theorem.m2r.QCM}
Given a smooth mapping $F:\R^n\rightarrow\R^m$, a closed-graph mapping $C:\R^n\rightrightarrows\R^m$, a point $(\bar u,\bar y)\in\gph\,C$, and a unit direction $w\in \R^n$ with $DC(\bar u|\bar y)(w)\neq\varnothing$, define the new set-valued mapping 
$\widetilde{\varXi}:\gph\,C \times\R^m\rightrightarrows\R^n$ by
\begin{align*}
\widetilde{\varXi}((u,y),z):=\nabla F(u) z+D^*C(u|y)(z).
\end{align*}
Then the following are equivalent:\vspace*{-0.05in}
\begin{enumerate}
\item[\bf(a)] $F+C$ is metrically $2$-regular around $(\bar u,F(\bar u)+\bar y)$ relative to $w$.

\item[\bf(b)] It holds that
\begin{align*}
\big[0\in \widetilde{\varXi}((u,y),z),\quad z\neq 0\big]\Longrightarrow \inf_{\eta\in DC(\bar u|\bar y)(w)}d\ell_{\widetilde{\varXi}}(\bar u,\bar y)\left(w,\eta\right)>0.
\end{align*}
\end{enumerate}
In particular, if $\ell_{\widetilde{\varXi}}(\bar u,\bar y)=0$, then the infimum of all $c>0$ for which $F+C$ is metrically 2-regular around $(\bar u,F(\bar u)+\bar y)$ relative to $w$ is equal to
\begin{align*}
\sup\limits_{\eta\in DC(\bar u|\bar y)(w)}\left(\frac{1}{ d\ell_{\widetilde{\varXi}}(\bar u,\bar y)\left(w,\eta\right)}\right).
\end{align*}
\end{theorem}
{\bf Proof}.
Combine Lemma~\ref{lemma.coderivative.sum} with Theorem~\ref{theorem.m2r.criterion}.\qed \vspace*{0.05in}

The obtained characterizations of metric 2-regularity are based on the coderivative criterion for metric regularity married to the subderivative calculation for the corresponding LSV function. In general, computing this subderivative is challenging, but we can employ Theorems~\ref{theorem.subderivative.singular.I} and \ref{theorem.subderivative.singular.II} to derive constructive sufficient conditions for metric 2-regularity. The lemma below, taken from \cite[Example~4.21]{JDiss24}, provides an application of Theorem~\ref{theorem.subderivative.singular.I} that allows us to ensure the metric 2-regularity of set-valued mappings represented as sums of smooth and constant {\em polyhedral} mappings.\vspace*{-0.05in}

\begin{lemma}\label{proposition.m2r.Constraint.Polyhedral}
Given a nonempty polyhedral set $C_0\subset\R^m$, a continuously differentiable mapping $F:\R^n\rightarrow\R^m$, and a point $(\bar u,\bar y)\in\R^n\times C_0$, suppose that $\nabla F$ is semidifferentiable at $\bar u$ for a unit direction $w\in\R^n$. Then  the following condition is sufficient for the metric $2$-regularity of $F(\cdot)+C_0$ around $(\bar u,F(\bar u)+\bar y)$ relative to $w$:
\begin{align}\label{f.m2r.constraint.map.polyhedral}
\big[z\in \left(-N_{C_0}(\bar y)\right)\bigcap\ker\nabla F(\bar u),\quad(\nabla F)'(\bar u;w) z\in \nabla F(\bar u) N_{C_0}(\bar y)\big]\Longrightarrow z=0.
\end{align}
\end{lemma}

The next theorem employs Lemma~\ref{proposition.m2r.Constraint.Polyhedral} to obtain a stronger result.\vspace*{-0.07in}

\begin{theorem}\label{theorem.m2r.polyhedral}
Given a polyhedral mapping $C:\R^n\rightrightarrows\R^m$, a continuously differentiable mapping $F:\R^n\rightarrow\R^m$, and a point $(\bar u,\bar y)\in \gph\,C$, assume that $\nabla F$ is semidifferentiable at $\bar u$ for a unit direction $w\in \dom(DC(\bar u|\bar y))$. Then the following condition is sufficient for the metric $2$-regularity of $F+C$ around $(\bar u,F(\bar u)+\bar y)$ relative to $w$:
\begin{align}\label{f.m2r.map.polyhedral}
\left.
\begin{array}{l}
0\in\nabla F(\bar u) z + D^*C(\bar u|\bar y)(z),\\ 
0\in(\nabla F)'(\bar u;w)z+\nabla F(\bar u)\nu+ D^*C(\bar u|\bar y)(\nu)\end{array}\right\}\Longrightarrow z=0.
\end{align}
\end{theorem}
{\bf Proof}. Denote $\varSigma:=F+C$ and $\bar v:=F(\bar u)+\bar y$. Nothing needs to be proved when $\varSigma$ is metrically regular around $(\bar u,\bar v)$, and thus it remains to consider the case where $\varSigma$ is not metrically regular around this point.

Denote further $G(u):=(-u,F(u))$, $\mathcal{C}(u):=G(u)+\gph\,C$ and observe by Theorem~\ref{theorem.metric.regularity.QCM} that the number $\textnormal{Reg}(u,F(u)+y;\varSigma)$ cannot be smaller than $\textnormal{Reg}(u,(0,F(u)+y);\mathcal{C})$ for any $(u,y)\in\gph\,C$. Moreover, it follows from the definition of the LSV function that the latter relation holds trivially for any $(u,y)\notin\gph\,C$. Thus we get
\begin{align}\label{p6.1}
\textnormal{Reg}(u,F(u)+y;\varSigma)\geq\textnormal{Reg}(u,(0,F(u)+y);\mathcal{C})\;\mbox{ for all }\;(u,y)\in \R^n\times\R^m.
\end{align}
Remembering that $\varSigma$ is not metrically regular around $(\bar u,\bar v)$ brings us by Theorem~\ref{theorem.metric-regularity} and Theorem~\ref{theorem.metric.regularity.QCM} to the equalities
\begin{align}\label{p6.2}
\textnormal{Reg}(\bar u,\bar v;\varSigma)=\textnormal{Reg}(\bar u,(0,\bar v);\mathcal{C})=0.
\end{align}
This allows us to deduce from the estimate in \eqref{p6.1} that the condition
\begin{align}\label{p6.3}
d\textnormal{Reg}(\bar u,\bar v;\varSigma)(w,\mu)\geq d\textnormal{Reg}(\bar u,(0,\bar v);\mathcal{C})\left(w,(0,\mu)\right)\;\mbox{ for all }\;\mu\in \R^m
\end{align}
is satisfied. Since $C$ is polyhedral, the mapping $\mathcal{C}(\cdot)=G(\cdot)+\gph\,C$ fits the framework of Lemma~\ref{proposition.m2r.Constraint.Polyhedral}, and hence we get from the above that the following condition is sufficient for the metric 2-regularity of $\mathcal{C}$ around $(\bar u,(0,\bar v))$ relative to $w$:
\begin{align*}
\left.
\begin{array}{r l}
\zeta & \in \left( -N_{\gph\,C}(\bar u,\bar y)\right)\bigcap \ker \nabla G(\bar u),\\
(\nabla G)'(\bar u;w) \zeta & \in \nabla G(\bar u) N_{\gph\,C}(\bar u,\bar y)
\end{array}\right\}\Longrightarrow\zeta=0.
\end{align*}
By the construction of $G$, we see that the obtained condition corresponds to
\begin{align*}
\left.
\begin{array}{r l}
(-\nabla F(\bar u) z,-z) & \in N_{\gph\,C}(\bar u,\bar y,\\
(-(\nabla F)'(\bar u;w) z-\nabla F(\bar u)\nu,-\nu) & \in N_{\gph\,C}(\bar u,\bar y)
\end{array}\right\}
\Longrightarrow z=0.
\end{align*}
The coderivative definition \eqref{cod} tells us that the latter is nothing else than \eqref{f.m2r.map.polyhedral}, and so we deduce from Lemma~\ref{proposition.m2r.Constraint.Polyhedral} that $\mathcal{C}$ is metrically 2-regular around $(\bar u,(0,\bar v))$ relative to $w$. It follows from \eqref{p6.2}, \eqref{p6.3}, and Theorem~\ref{theorem.m2r.criterion} that
\begin{align}\label{p6.4}
\left(\inf_{(0,\mu)\in D\mathcal{C}(\bar u|(0,\bar v))(w)}d\textnormal{Reg}(\bar u,\bar v;\varSigma)(w,\mu)\right)>0.
\end{align}
At the same time, Lemma~\ref{lemma.coderivative.sum} and \cite[Example~3.19]{JDiss24} yield
\begin{align*}
D\mathcal{C}(\bar u|(0,\bar v))(w)=\left(\begin{array}{c}
-w\\
F'(u^*)w
\end{array}\right)+T_{\gph\,C}(\bar u,\bar y),
\end{align*}
i.e., we have $(0,\mu)\in D\mathcal{C}(\bar u|(0,\bar v))(w)$ if and only if
\begin{align*}
0=-w+\zeta,\quad \mu=F'(u^*) w+\lambda,\quad \lambda\in DC(\bar u|\bar y)(\zeta).
\end{align*}
Thus condition \eqref{p6.4} can be rewritten as
\begin{align*}
\left(\inf_{\mu\in D\varSigma(\bar u|\bar v)(w)}d\textnormal{Reg}(\bar u,\bar v;\varSigma)(w,\mu)\right)>0,
\end{align*}
which allows us to deduce from \eqref{p6.2} that assertion (c) of Theorem~\ref{theorem.m2r.criterion} is fulfilled. Since all the assumptions of that theorem are satisfied, we conclude that $\varSigma=F+C$ is metrically 2-regular around $(\bar u,\bar v=F(\bar u+\bar y))$ relative to $w$.\qed\vspace*{0.05in}

Let us specify the obtained theorem for the case where the polyhedral mapping $C(\cdot)$ is the {\em indicator mapping} \eqref{indic} associated with a polyhedral set.\vspace*{-0.05in}

\begin{corollary}\label{corollary.m2r.Indicator.Polyhedral} Given  a nonempty polyhedral set $C_0\subset\R^n$, a continuously differentiable mapping $F:\R^n\rightarrow\R^m$, and a point $u^*\in C_0$, assume that $\nabla F$ is semidifferentiable at $\bar u$ for a unit direction $w\in T_{C_0}(\bar u)$. Then the following condition is sufficient for the metric 2-regularity of $F+\varDelta_{C_0}$ around $(\bar u,F(\bar u))$ relative to $w$:
\begin{align}\label{f.m2r.indicator.map.polyhedral.A2}
\big[\nabla F(\bar u) z\in N_{C_0}(\bar u),\;(\nabla F)'(\bar u;w) z\in \textnormal{rge}\,\nabla F(\bar u)+N_{C_0}(\bar u)\big]\Longrightarrow z=0.
\end{align}
\end{corollary}
\begin{proof}
This is a direct consequence of Lemma~\ref{prop.coderivative.calculus} and Theorem~\ref{theorem.m2r.polyhedral}.\qed
\end{proof}\vspace*{-0.15in}

\begin{remark}\label{remark.m2r.2regwrtS} It is shown  in \cite[Remark~4.22]{JDiss24} that if $F$ is twice differentiable, $C_0\subset\R^m$ is polyhedral, and $F(\bar u)=\bar y\in C_0$, then condition \eqref{f.m2r.constraint.map.polyhedral} agrees with the notion of 2-\textit{regularity of $F$ around $\bar u$ with respect to the set $C_0$} for the unit direction $w$ formulated in \cite[Definition~1]{AAI05} as follows:
\begin{align}\label{f.2r.wrts}
\textnormal{rge}\,F'(\bar u)+F''(\bar u)\left[w,F'(\bar u)^{-1}T_{C_0}(\bar y)\right]-T_{C_0}(\bar y)=\R^m,
\end{align}
where $F'(\bar u)=\nabla F(\bar u)^\top$, and where $F''(\bar u)[w,\cdot]:\R^n\rightarrow\R^m$ is the linear mapping associated with the second derivative of $F$ at $\bar u$ relative to $w$.  Using \cite[Theorem~13.2]{Rocka}, we can write this down in our notation as
\begin{align*}
(\nabla F)' (\bar{u};w)^\top\left(\nabla F(\bar u)^{-\top}T_{C_0}(\bar y)\right)=F''(\bar u)\left[ w,F'(\bar u)^{-1}T_{C_0}(\bar y)\right].
\end{align*}\qed
\end{remark}\vspace*{-0.1in}

\begin{remark}\label{remark.m2r.2regwrtS.2} It is argued in \cite[Remark~4.25]{JDiss24} that if $F$ is twice differentiable, $C_0\subset\R^n$ is polyhedral, and $\bar u\in C_0$, then  our condition \eqref{f.m2r.indicator.map.polyhedral.A2} is generally stronger than the condition in \cite[Theorem~5.1]{AruIz20} for the 2-{\em regularity of $F$ around $\bar u$ with respect to $C_0$} in the unit direction $w$ formulated as
\begin{align}\label{f.2r.wrts.II}
F'(\bar u)\mathcal{L}+F''(\bar u)\left[ w,\ker F'(\bar u)\cap T_{C_0}(\bar u)\right]=\R^m,
\end{align}
where $F'(\bar u)=\nabla F(\bar u)^\top$, $\mathcal{L}$ is the linear hull (span) of the set $C_0$, and $F''(\bar u)[w,\cdot]:\R^n\rightarrow\R^m$ is the mapping from Remark~\ref{remark.m2r.2regwrtS}. It would be interesting to check whether \eqref{f.2r.wrts.II} is sufficient for the metric 2-regularity of $F+\Delta_{C_0}$, since this could improve the outcome of Corollary~\ref{corollary.m2r.Indicator.Polyhedral}. At the same time, conditions  \eqref{f.m2r.indicator.map.polyhedral.A2} and \eqref{f.2r.wrts.II} agree if
\begin{align*}
\nabla F(\bar u)^{-1}\mathcal{L}^\bot= \nabla F(\bar u)^{-1} N_{C_0}(\bar u),
\end{align*}
which holds, in particular, if either $F'(\bar u)=0$, or $\bar u$ is in the relative interior of $C_0$.\qed
\end{remark}\vspace*{-0.06in}

The latter remark says that \eqref{f.m2r.indicator.map.polyhedral.A2} is sufficient for the 2-regularity condition \eqref{f.2r.wrts.II} in the case of twice differentiable mappings. The following example shows that the converse implication is not generally true. Moreover, this example tells us that \eqref{f.m2r.indicator.map.polyhedral.A2} is only sufficient but not necessary for the metric 2-regularity under consideration.\vspace*{-0.05in}

\begin{example}\label{ex.2reg.1} Let $F(u):=u$, $C_0:=\R_+$, $\bar u:=0$, and $w:=1$. We have $\mathcal{L}=\R$ and $N_{C_0}(\bar u)=\R_{-}$, so \eqref{f.2r.wrts.II} holds, while \eqref{f.m2r.indicator.map.polyhedral.A2} does not. To see that \eqref{f.m2r.indicator.map.polyhedral.A2} is not necessary for the  metric 2-regularity, observe that 
\begin{align*}
\textnormal{Reg}(u,F(u);F+\Delta_{C_0})=1\;\mbox{ for all }\;u>0,
\end{align*}
so Theorem~\ref{theorem.m2r.criterion} gives metric 2-regularity of $F+\Delta_{C_0}$ around $(\bar u,F(\bar u))$ relative to $w$.\qed
\end{example}\vspace*{-0.05in}

Lemma~\ref{proposition.m2r.Constraint.Polyhedral}, Theorem~\ref{theorem.m2r.polyhedral}, and Corollary~\ref{corollary.m2r.Indicator.Polyhedral} can be regarded as applications of Theorem~\ref{theorem.subderivative.singular.I}. In the rest of this section, we consider applications of Theorem~\ref{theorem.subderivative.singular.II}. The next result is taken from \cite[Example~4.27]{JDiss24}. \vspace*{-0.05in}

\begin{proposition}\label{proposition.m2r.Constraint.NonPolyhedral}
Given a closed set $C_0\subset\R^m$, a continuously differentiable mapping $F:\R^n\rightarrow\R^m$, and a point $(\bar u,\bar y)\in\R^n\times C_0$, assume that $\nabla F$ is semidifferentiable at $\bar u$ for a unit direction $w\in\R^n$. 
Then the following condition is sufficient for the metric 2-regularity of $F(\cdot)+C_0$ around $(\bar u,F(\bar u)+\bar y)$ relative to $w$:
\begin{align}\label{f.m2r.constraint.map.nonpolyhedral}
\Big[z\in \left(-N_{C_0}(\bar y)\right)\bigcap \ker \nabla F(\bar u),\;(\nabla F)'(\bar u;w) z\in\textnormal{rge}\nabla F(\bar u)\Big]\Longrightarrow z=0.
\end{align}
\end{proposition}

The latter proposition generalizes Lemma~\ref{proposition.m2r.Constraint.Polyhedral} in the sense that $C_0$ is not necessarily polyhedral. At the same time, \eqref{f.m2r.constraint.map.nonpolyhedral} may be stronger than \eqref{f.m2r.constraint.map.polyhedral} in Lemma~\ref{proposition.m2r.Constraint.Polyhedral}, and Example~\ref{example.nom2r.nonpoly} below demonstrates that \eqref{f.m2r.constraint.map.nonpolyhedral} cannot be simply replaced by the weaker condition \eqref{f.m2r.constraint.map.polyhedral} without destroying the conclusion of 
Proposition~\ref{proposition.m2r.Constraint.NonPolyhedral}.\vspace*{0.03in}

Now we establish a novel sufficient condition for metric 2-regularity, which can be viewed as a nonpolyhedral counterpart of Corollary~\ref{corollary.m2r.Indicator.Polyhedral}. \vspace*{-0.05in}

\begin{theorem}\label{theorem.m2r.Indicator.nonPolyhedral}
Given a closed set $C_0\subset\R^n$, a continuously differentiable mapping $F:\R^n\rightarrow\R^m$, and a point $\bar u\in C_0$, suppose that
\begin{align}\label{f.mr.pre}
\textnormal{rge}\nabla F(\bar u)\bigcap N_{C_0}(\bar u)=\{ 0\}.
\end{align}
If $\nabla F$ is semidifferentiable at $\bar u$ for a unit direction $w\in T_{C_0}(\bar u)$, 
then condition \eqref{f.m2r.indicator.map.polyhedral.A2} is sufficient for the metric $2$-regularity of $F+\varDelta_{C_0}$ around  $(\bar u,F(\bar u))$ relative to $w$.
\end{theorem}\vspace*{-0.05in}
{\bf Proof}. Consider the  closed set
\begin{align}\label{f.subderivative.singular.nonpolyhedral.Z0}
Z_0:=\left\{ z\in\R^m\mid F'(u^*)^\top z\in N_{C_0}(\bar u)\right\}
\end{align}
and observe that if $Z_0=\{ 0\}$, then Lemmas~\ref{lemma.coderivative.sum} and \ref{prop.coderivative.calculus} yield the metric regularity, and hence metric 2-regularity relative to $w$, of the mapping $F+\varDelta_{C_0}$ around $(\bar u,F(\bar u))$. Assuming from now on that 
\begin{align}\label{p8.1}
Z_0\setminus\{ 0\}\neq\varnothing
\end{align}
for the mapping $\varDelta_{C_0}$, we verify conditions (i)--(iii) of Theorem~\ref{theorem.subderivative.singular.I} with $\mathcal{D}:=\gph\,\varDelta_{C_0}=C_0\times \{ 0\}$, $\mathcal{A}(\xi):=\mathcal{A}(u,y):=\nabla F(u)$, $\varGamma(\xi,z):=\varGamma((u,y),z):= D^*\varDelta_{C_0}(u|0)(z)=N_{C_0}(u)$, and $\omega:=(w,\eta)$ as $\eta \in D\varDelta_{C_0}(\bar u|0)(w)=\varDelta_{T_{C_0}(u^*)}(w)=\{ 0\}$, where the last equalities follow by direct computations. Indeed, we know that $\varGamma$ is outer semicontinuous, and so condition (ii) of Theorem~\ref{theorem.subderivative.singular.I} holds. As stated above, condition \eqref{p8.1} means that $F+\varDelta_{C_0}$ is not metrically regular around $(\bar u,F(\bar u))$. Thus the equality
$\textnormal{Reg}(\bar u,F(\bar u);F+\varDelta_{C_0})=0$
follows by Theorem~\ref{theorem.metric-regularity} implying in turn that condition (i) of Theorem~\ref{theorem.subderivative.singular.I} is fulfilled. Condition (iii) is a consequence of the semidifferentiability of $\nabla F$ at $\bar u$ for $w$. The mapping $\varGamma$ satisfies the assumptions of Proposition~\ref{prop1.suff.LSV.subderivative.II}, which ensures that condition (iv') of Theorem~\ref{theorem.subderivative.singular.II} is satisfied. Condition (v') of that theorem is a direct consequence of \eqref{f.mr.pre}, which is guaranteed by the imposed assumptions. Therefore, all the conditions of Theorem~\ref{theorem.subderivative.singular.II} are fulfilled, and we get
\begin{align*}
d\ell_{\widetilde{\varXi}}(u^*,0)(w,0) \geq
\min_{\substack{z\in Z_0\\\Vert z\Vert=1}}\left( \di\left[0, (\nabla F)'(\bar u;w) z+\textnormal{rge}\nabla F(\bar u) -N_{C_0}(\bar u)\right]\right),
\end{align*}
where $\widetilde{\varXi}$ is defined as in Theorem~\ref{theorem.m2r.QCM}. With \eqref{f.m2r.indicator.map.polyhedral.A2}, \eqref{p8.1}, and $D\varDelta_{C_0}(\bar u|0)(w)=\{ 0\}$ in mind, 
the metric 2-regularity of $F+\varDelta_{C_0}$ is finally deduced from Theorem~\ref{theorem.m2r.QCM}~(b).\qed \vspace*{0.05in}

Note that condition \eqref{f.mr.pre} does not necessitate the metric regularity of $F+\varDelta_{C_0}$ around $(\bar u,F(\bar u))$; see \cite[Section~3]{FIJ24} for more discussions. Let us illustrate by the example below that \eqref{f.m2r.indicator.map.polyhedral.A2} alone is not enough to guarantee the metric 2-regularity. This means that \eqref{f.mr.pre} cannot be removed from the assumptions without destroying the conclusion of Theorem~\ref{theorem.m2r.Indicator.nonPolyhedral}, and this is true even if the set $C_0$ is a smooth manifold.\vspace*{-0.07in}

\begin{example}\label{ex.indicator.remove.separation}
Consider the case where $C_0:=h^{-1}(0)$ for $h(u)=h(u_1,u_2,u_3):=(1/2)(u_2+u_3^2)$, $\bar u:=(0,0,0)$, and the mapping $F$ is given by
\begin{align*}
F(u):=F(u_1,u_2,u_3):=\left(\begin{array}{c}
u_1^2+u_2+u_3^2\\
u_1
\end{array}
\right).
\end{align*}
The set $C_0$ is a smooth manifold, and it is easy to check that
\begin{align*}
T_{C_0}(u) & =\ker\,h'(u)=\left\{ (w_1,w_2,w_3)\,\left|\, \frac{1}{2}w_2+u_3w_3=0\right.\right\},\\
N_{C_0}(u) & =\textnormal{rge}\,\nabla h(u)=\textnormal{span}\left\{ \left( 0, \frac{1}{2},u_3\right)\right\}\;\mbox{ for all }\;u\in C_0.
\end{align*}
Furthermore, we have for $w:=(0,0,1)$ that
\begin{align*}
\nabla F(u) & =\left(\begin{array}{c c}
2u_1 & 1\\
1 & 0\\
2u_3 & 0
\end{array}\right),\qquad (\nabla F)'(\bar u; w)=\left(\begin{array}{c c}
0 & 0\\
0 & 0\\
2 & 0
\end{array}
\right)
\end{align*}
whenever $u\in C_0$. Observe that $w\in T_{C_0}(\bar u)$ and condition \eqref{f.m2r.indicator.map.polyhedral.A2} is satisfied while \eqref{f.mr.pre} fails. To show by using Theorem~\ref{theorem.m2r.criterion} that $\varSigma:=F+\varDelta_{C_0}$ cannot be metrically 2-regular around $(\bar u,F(\bar u))$ relative to $w$, pick any $t\searrow 0$, put  $w(t):=(0,-t,1)$, find $w(t)\rightarrow w$ as $t\searrow 0$,  and denote $u(t):=\bar u+tw(t)\in C_0$. We have $v(t):=-(0,1,2t)\in N_{C_0}(u(t))$ for all $t$ and then get with $\bar z:=(1,0)$ that
\begin{align*}
\Reg(u(t),F(u(t));\varSigma) & =\inf_{\Vert z\Vert=1}\di\left[0,\nabla F(u(t))z+N_{C_0}(u(t))\right]\\
~ &\leq \di\left[0,\nabla F(u(t))\bar{z}+N_{C_0}(u(t))\right]\\
~ & \leq\left\Vert(0,1+v_2(t),2t+v_3(t))\right\Vert\equiv 0.
\end{align*}
This tells us that Theorem~\ref{theorem.m2r.criterion} (c) is violated, and thus $\varSigma=F+\varDelta_{C_0}$ cannot be metrically 2-regular around the point in question relative to $w$.\qed
\end{example}\vspace*{-0.05in}

The next example illustrates, again for smooth manifolds, that condition \eqref{f.m2r.constraint.map.nonpolyhedral} of Proposition~\ref{proposition.m2r.Constraint.NonPolyhedral} cannot be replaced by \eqref{f.m2r.constraint.map.polyhedral} without destroying the conclusion.\vspace*{-0.05in}

\begin{example}\label{example.nom2r.nonpoly} In the setting of Example~\ref{ex.indicator.remove.separation}, let $G(u):=(-u,F(u))$ and $D_0:=\gph\,\varDelta_{C_0}=C_0\times\{ (0,0)\}$. It follows for $\bar u=(0,0,0)$ that
\begin{align*}
N_{D_0}(\bar u,(0,0))=\{ 0\}\times\R\times\{ 0\}\times\R^2.
\end{align*}
Furthermore, for $w=(0,0,1)$ we get
\begin{align*}
\nabla G(\bar u)=\left(\begin{array}{c c c c c}
-1 & 0 & 0& 0 & 1\\
0 & -1 & 0 & 1 & 0\\
0 & 0 & -1 & 0 & 0
\end{array}\right),\quad (\nabla G)'(\bar u; w)=\left(\begin{array}{c c c c c}
-1 & 0 & 0& 0 & 1\\
0 & -1 & 0 & 1 & 0\\
0 & 0 & -1 & 2 & 0
\end{array}\right)
\end{align*}
This allows us to observe that \eqref{f.m2r.constraint.map.polyhedral} is satisfied while \eqref{f.m2r.constraint.map.nonpolyhedral} fails. Considering for $t\searrow 0$ the vector $u(t)$ in Example~\ref{ex.indicator.remove.separation},  we deduce from
Theorem~\ref{theorem.metric.regularity.QCM} with $\varSigma=F+\varDelta_{C_0}$ and $\mathcal{C}(\cdot)=G(\cdot)+D_0$ that
$0\equiv\Reg(u(t),F(u(t));\varSigma)\geq \Reg(u(t),(0,F(u(t)));\mathcal{C})\geq 0$. This tells us by Theorem~\ref{theorem.m2r.criterion} (c) that $\mathcal{C}$ is not metrically 2-regular around the given point relative to $w$. The only possible reason for this is violating \eqref{f.m2r.constraint.map.nonpolyhedral}.\qed
\end{example}\vspace*{-0.07in}

The final result of this section concerns {\em product mappings} between product spaces.\vspace*{-0.07in}

\begin{theorem}\label{theorem.m2r.productMap}
Given a continuously differentiable mapping $G:\R^k\times\R^l\rightarrow\R^{p+q}$ and outer semicontinuous mappings $R:\R^k\rightrightarrows\R^p$ and $T:\R^l\rightrightarrows\R^q$, consider the product mapping $\mathcal{P}:\R^k\times\R^l\rightrightarrows\R^m$ from \eqref{f.product.map}. Let a point $((\bar x,\bar \sigma),(\bar y,\bar\nu))\in\gph\,\mathcal{P}$ and a unit direction $(w,\mu)\in \R^k\times\R^l$ with $D\mathcal{P}\left((\bar x,\bar\sigma)|(\bar y,\bar\nu)\right)(w,\mu)\neq\varnothing$ satisfy the following:\vspace*{-0.05in}
\begin{enumerate}
\item[\bf(a)] $\nabla G$ is semidifferentiable at $(\bar x,\bar\sigma)$ for $(w,\mu)$.
\item[\bf(b)] $T$ is Lipschitz-like around $(\bar\sigma,\bar\nu)$.
\item[\bf(c)] $T$ is metrically regular around $(\bar\sigma,\bar\nu)$.
\item[\bf(d)] We have the conditions
\begin{align*}
\rge(D^* R(\bar x,\bar y))=\{ 0\},\qquad \textnormal{rge}\left(\nabla_\sigma G(\bar x,\bar\sigma)\right)\bigcap\rge(D^*T(\bar\sigma|\bar\nu))=\{ 0\}. 
\end{align*}
\item[\bf(e)] The pair $z=(0,0)$ is the only solution $z=(\theta,\zeta)\in\R^{p+q}$ to the system
\begin{align*}			
0 & =\nabla_x G(\bar x,\bar\sigma)z,\\
0 & \in \nabla_\sigma G(\bar x,\bar\sigma) z+D^*T(\bar\sigma|\bar\nu)(\zeta),\\
0 & = (\nabla_x G)'\left((\bar x,\bar\sigma);(w,\mu)\right) z +\nabla_x G(\bar x,\bar\sigma)\beta,\\
0 & \in (\nabla_\sigma G)'\left((\bar x,\bar\sigma);(w,\mu)\right)z + \nabla_\sigma G(\bar x,\bar\sigma)\beta+\rge(D^*T(\bar\sigma|\bar \nu)).
\end{align*}
\end{enumerate}
Then the summation mapping $G+\mathcal{P}$ is metrically 2-regular around $((\bar x,\bar\sigma),G(\bar x,\bar\sigma)+(\bar y,\bar\nu))$ relative to $(w,\mu)$.
\end{theorem}\vspace*{-0.05in}
{\bf Proof}. First we are going to check that all the conditions of Theorem~\ref{theorem.subderivative.singular.II} are satisfied under the assumptions imposed now. To that end, define
\begin{align*}
Z_0:=\left\{ z\in\R^m\;\left|\;0\in \nabla G(\bar x,\bar \sigma) z+D^*\mathcal{P}((\bar x,\bar\sigma)|(\bar y,\bar\nu))(z)\right.\right\},
\end{align*}
and recall that the mapping $G+\mathcal{P}$ is metrically regular around $\big((\bar x,\bar\sigma),G(\bar x,\bar\sigma)+(\bar y,\bar\nu)\big)$ if and only if $Z_0=\{ 0\}$. There is nothing to prove in this case, so we suppose that $Z_0\neq\varnothing$. Pick any $z=(\theta,\zeta)\in\R^{p+q}$ and deduce from Lemma~\ref{lemma.coderivative.product} that
\begin{align*}
D^*\mathcal{P}((\bar x,\bar\sigma)|(\bar y,\bar\nu))(z)=D^*R(\bar x|\bar y)(\theta)\times D^*T(\bar\sigma|\bar \nu)(\zeta).
\end{align*}
According to the first equality in (d), we have
\begin{align}\label{p9.1}
D^*\mathcal{P}((\bar x,\bar\sigma)|(\bar y,\bar\nu))(z)=\{ 0\}\times D^*T(\bar y|\bar\nu)(\zeta)
\end{align}
for all such vectors $z$. This allows us to write
\begin{align}\label{p9.12}
Z_0=\left\{ z=(\theta,\zeta)\;\left|\;0=\nabla_x G(\bar x,\bar\sigma)z,\;
0\in \nabla_\sigma G(\bar x,\bar\sigma) z+D^*T(\bar\sigma|\bar\nu)(\zeta)\right.\right\}.
\end{align}
The fulfillment of conditions (i)--(iii) of Theorem~\ref{theorem.subderivative.singular.II} can be checked by our standard arguments already employed in the proof of 
Theorem~\ref{theorem.m2r.polyhedral}. 
Theorem~\ref{theorem.metric-regularity} tells us that conditions (b) and (c) are equivalent to
\begin{align}\label{p9.2}
D^*T(\bar y|\bar\nu)(0)=\{ 0\},\qquad D^*T(\bar y|\bar\nu)^{-1}(0)=\{0\}.
\end{align}
Therefore, using \eqref{p9.1} and Proposition~\ref{prop2.suff.LSV.subderivative.II} yields the fulfillment of condition (iv') of Theorem~\ref{theorem.subderivative.singular.II}. Condition (v') of that theorem transforms in our case into 
\begin{align}\label{p9.3}
\textnormal{rge}\nabla G(\bar x,\bar\sigma)\bigcap\varTheta\subset\{ 0\},
\end{align}
where the closed cone $\varTheta$ is defined by
\begin{align*}
\varTheta:=\cl\left(\cone\left( \bigcup_{\substack{\Vert z\Vert=1\\
z\in\dom(D^*\mathcal{P}((\bar x,\bar\sigma)|(\bar y,\bar\nu)))}} D^*\mathcal{P}((\bar x,\bar\sigma)|(\bar y,\bar\nu))(z)\right)\right).
\end{align*}
It follows from \eqref{p9.1}, \eqref{p9.2}, and Lemma~\ref{lemma.closed.domain.range} that 
\begin{align}\label{p9.4}
\varTheta=\{ 0\}\times \rge\left(D^*T(\bar\sigma|\bar\nu)\right),
\end{align}
and so \eqref{p9.3} transforms into
$\textnormal{rge}\nabla_\sigma G(\bar x,\bar\sigma)\bigcap\rge\left(D^*T(\bar\sigma|\bar\nu)\right)\subset\{ 0\}$,
which is guaranteed by the second equality in condition (d) of the theorem. Therefore, all the conditions of Theorem~\ref{theorem.subderivative.singular.II} are fulfilled, and we get for all $\eta\in\R^m$ that
\begin{align*}
~ & d\ell_{\widetilde{\varXi}}\left((\bar x,\bar\sigma),(\bar y,\bar\nu)\right)((w,\mu),\eta) \\
\geq & 
\min_{\substack{z\in Z_0\\\Vert z\Vert=1}}\left( \di\left[0, (\nabla G)'((\bar x,\bar\sigma);(w,\mu)) z+\varTheta+\textnormal{rge}\nabla G(\bar x,\bar\sigma)\right]\right),
\end{align*}
where the mapping $\widetilde{\varXi}$ is defined as in Theorem~\ref{theorem.m2r.QCM}. With \eqref{p9.12} and \eqref{p9.4} in mind, it follows from condition (e) of this theorem that
\begin{align*}
\left(\inf_{\eta\in D\mathcal{P}((\bar x,\bar\sigma)|(\bar y,\bar\nu))}d\ell_{\widetilde{\varXi}}\left((\bar x,\bar\sigma),(\bar y,\bar\nu)\right)((w,\mu),\eta)\right)>0.
\end{align*}
Hence, we get the claimed assertion by employing Theorem~\ref{theorem.m2r.QCM}.\qed\vspace*{-0.2in}

\section{Gfrerer Regularity}\label{sec.Gfr}\vspace*{-0.1in}

In this section, we consider Gfrerer's concept of \textit{metric pseudo-regularity of order 2} introduced in \cite[Definition~1.1]{Gf13}. In that paper and in \cite{BM25}, Gfrerer's concept is used to develop tractable necessary optimality conditions for constrained optimization problems, while in \cite{JDiss24,J24} it is used in conjunction with second-order contingent derivatives \cite{DuFlo18,J24a}  
to calculate tangents to a level set of a mapping in the absence of standard first-order regularity conditions.
We will outline below further possible applications of Gfrerer's concept and discuss relations with metric 2-regularity.\vspace*{-0.05in}

\begin{definition}\label{definition.gfrerer.property}
A mapping $S:\R^n\rightrightarrows\R^m$ is said to be \textit{Gfrerer regular} around $(\bar u.\bar y)\in\gph\,S$ relative to a direction $(w,\eta)\in\R^n\times\R^m$ if  there exist $\varepsilon,\delta,c>0$ such that \begin{align*}
\di\left[u,S^{-1}(y)\right]\cdot \Vert u-\bar u \Vert\leq c\cdot \di[y,S(u)]
\end{align*}
for all $(u,y)\in K_{\varepsilon,\delta}\left((\bar u,\bar y);(w,\eta)\right)$ with $\di[y,S(u)]\leq \delta\Vert u-\bar u\Vert^2$. 
\end{definition}\vspace*{-0.05in}

The discussion after \cite[Definition~1]{Gf13} shows that $S$ is Gfrerer regular around $(\bar u,\bar y)$ relative to $(w,\eta)$ whenever $\eta\notin DS(\bar u|\bar y)(w)$. The next result concerns Gfrerer regularity in the {\em opposite case}, which is relative to directions $(w,\eta)\in T_{\gph\,S}(\bar u,\bar y)$. Again, the {\em LSV function of the coderivative} is our driving force here.\vspace*{-0.05in}

\begin{theorem}\label{theorem.gamma.regularity.gfrerer}
Let $S:\R^n\rightrightarrows\R^m$ be a set-valued mapping whose graph is closed around a point $(\bar u,\bar y)\in\gph\,S$, and let $(w,\eta)\in T_{\gph\,S}(\bar u,\bar y)$ be a given direction. Then the following assertions are equivalent:\vspace*{-0.05in}
\begin{enumerate}
\item[\bf(a)] $S$ is Gfrerer regular around $(\bar u.\bar y)$ relative to $(w,\eta)$.

\item[\bf(b)] There are numbers $\rho,\kappa>0$ such that assertions {\rm(a)--(d)} of Lemma~\ref{lemma.gamma-regularity} are satisfied for the mapping $\varPsi:\R^n\rightrightarrows\R^m$ defined by
\begin{align*}
\gph\,\varPsi=\gph\,S\cap K_{\rho,\kappa}((\bar u,\bar y);(w,\eta)).
\end{align*}
\item[\bf(c)] We have the implication
\begin{align}\label{f.gfrerer.criterion}
\textnormal{Reg}(\bar u,\bar y;S)=0\Longrightarrow d\textnormal{Reg}(\bar u,\bar y;S)(w,\eta)>0.
\end{align}
\end{enumerate}
Furthermore, if $\textnormal{Reg}(\bar u,\bar y;S)=0$, then the infimum of all $c>0$ with which $S$ is Gfrerer regular around $(\bar u,\bar y)$ relative to $(w,\eta)$ is calculated by
$\big(d\textnormal{Reg}(\bar u,\bar y;S)(w,\eta)\big)^{-1}.$
\end{theorem}\vspace*{-0.05in}
{\bf Proof}. We can use \cite[Lemma~1]{Gf13} to verify implication (a)$\Rightarrow$(b) and \cite[Theorem~1]{Gf13} to get (c)$\Rightarrow$(a). Implication (b)$\Rightarrow$(c) essentially follows from Lemma~\ref{lemma.gamma-regularity}. More details and discussions can be found in \cite{JDiss24}. \qed\vspace*{0.05in}

The next result establishes a precise relationship between the notions of metric 2-regularity and Gfrerer regularity for set-valued mappings.\vspace*{-0.03in}

\begin{theorem}\label{theorem.m2r.gfrerer}
Let the graph of $S:\R^n\rightrightarrows\R^m$ be closed around $(\bar u,\bar y)\in\gph\,S$, and let $w\in\R^n$ be a unit direction. Then the following are equivalent:\vspace*{-0.05in}
\begin{enumerate}

\item[\bf(a)] $S$ is metrically 2-regular around $(\bar u,\bar y)$ relative to $w$.
\item[\bf(b)] $S$ is Gfrerer regular around $(\bar u,\bar y)$ relative to $(w,\eta)$ for all $\eta\in\R^m$, and we have
\begin{align*}
DS(\bar u|\bar y)(w)\neq\varnothing.
\end{align*}
\end{enumerate}
\end{theorem}
{\bf Proof}. We know that $S$ is Gfrerer regular around $ (\bar u,\bar y)$ relative to $(w,\eta)$ whenever $(w,\eta)\notin T_{\gph\,S}(\bar u,\bar y)$, i.e., $\eta\notin DS(\bar u|\bar y)(w)$. The definition of metric 2-regularity implies in turn that $DS(\bar u|\bar y)(w)\neq \varnothing$. Therefore, the claimed equivalence follows from comparing Theorem~\ref{theorem.m2r.criterion}(c) with Theorem~\ref{theorem.gamma.regularity.gfrerer}(c).\qed\vspace*{0.05in}

Metric 2-regularity can be viewed as a \textit{complete} version of Gfrerer regularity. Moreover, the obtained characterizations and sufficient conditions for metric 2-regu-larity guarantee the fulfillment of Gfrerer regularity.\vspace*{0.05in}

Finally, we discuss relationships between the notions and results developed above with those known in the literature.
\vspace*{-0.07in}

\begin{remark}\label{rem.2rg.classic}
It is argued in \cite[Example~6.4]{JDiss24} that for a continuously differentiable mapping $F:\R^n\rightarrow\R^m$ such that $\nabla F$ is semidifferentiable at $\bar u$ in a unit direction $w\in\R^n$, the following statements are equivalent:\vspace*{-0.07in}
\begin{enumerate}
\item[\bf(a)] $F$ is metrically 2-regular around $\bar u$ relative to $w$.
\item[\bf(b)] $F$ is Gfrerer regular around $\bar u$ relative to $(w,\eta)$ for all $\eta\in\R^m$.
\item[\bf(c)] We have the implication
\begin{align}\label{f.2-regularity.semi}\big[z\in\ker \nabla F(\bar u),\quad (\nabla F)'(\bar u;w) z\in\textnormal{rge}\nabla F(\bar u)\big] \Longrightarrow z=0.
\end{align}
\item[\bf(d)] We have the equality
\begin{align}\label{f.2reg.semi}
\textnormal{rge}F'(\bar u)+F''(\bar u;w)\ker F'(\bar u)=\R^m.
\end{align}
\end{enumerate}
If $\nabla F$ is locally Lipschitzian around $\bar u$  and semidifferentiable at this point for $w$, then \eqref{f.2reg.semi} agrees  with the \textit{2-regularity} condition in \cite{IzSol02} for differentiable mappings whose derivative is not necessarily differentiable. The latter 2-regularity condition is used in 
\cite{FIJ21,FIJ23B,Iz_CrSol_16} and elsewhere. This clearly
generalizes the classical 2-regularity notion in \cite{Av85,Tret84} designed for twice differentiable mappings. It has been first noted in \cite{Gf13} (see also \cite[Proposition~2]{GfOu16}) that for such mappings, Gfrerer regularity and classical 2-regularity are closely related properties.
Our observations extend this, in particular, to differentiable mappings whose derivative is merely semidifferentiable.\qed
\end{remark}\vspace*{-0.35in}

 \section{Applications to Parametric Constraint Systems}\label{sec.App}\vspace*{-0.1in}

 This section provides some applications of the new results obtained above for metric regularity, metric 2-regularity, and Gfrerer regularity of general set-valued mappings to particular classes of parametric problems governed by solution maps to {\em coupled constraint systems} of the type
 \begin{align}\label{f.CS}
 0=\varPhi(x,\sigma),\quad x\in\Omega,\quad 0\in T(\sigma),
 \end{align}
 where $\varPhi:\R^k\times\R^l\rightarrow\R^p$ is continuously differentiable, $\Omega\subset\R^k$ is closed, and $T:\R^l\rightrightarrows\R^q$ is outer semicontinuous. Systems of this type broadly appear in various theoretical settings and practical modeling; see, 
 e.g., \cite{BSh:PAOP:00,Rock_ImpF,FP03-1,ISo:NTMOVP:14,KK02,Mo06,Mo18,Mo24}. 

We rewrite the constraint system \eqref{f.CS} as a {\em generalized equation} of the form
\begin{align}\label{f.CS.reform}
 0\in G(x,\sigma)+\mathcal{P}(x,\sigma),
 \end{align}
 where $G:\R^k\times\R^l\rightarrow\R^{k+p+q}$ and $\mathcal{P}:\R^k\times\R^l\rightrightarrows\R^{k+p+q}$ are defined by
 \begin{align}
 G(x,\sigma) & :=\left(\begin{array}{c}-x\\\varPhi(x,\sigma)\\ 0\end{array}\right),\label{f.map.G}\\
 \mathcal{P}(x,\sigma) & :=\left(\begin{array}{c}\Omega\\ \{ 0\}\\ T(\sigma)\end{array}\right):=\left\{ \left.
 \left(\begin{array}{c}a\\0\\ b\end{array}\right)\right|\, a\in\Omega,\, b\in T(\sigma)
 \right\}.\label{f.map.P}
 \end{align}
 This setting fits the framework of the previous sections, 
 and we will now walk through our results to derive tractable conditions for the metric regularity, metric 2-regularity, and Gfrerer regularity of the mapping $G+\mathcal{P}$. 
 To this end, let us first compute derivatives of $G$ and $\mathcal{P}$ for a solution $(\bar x,\bar\sigma)$ to system \eqref{f.CS}. For $G$ in \eqref{f.map.G}, we have
\begin{align}\label{f.Part.der.G}
\nabla G(x,\sigma) =\left(\begin{array}{c c c}
 -\mathcal{I} & \nabla_x\varPhi(x,\sigma) & 0\\
 0 & \nabla_\sigma\varPhi(x,\sigma) & 0
 \end{array}\right),
\end{align}
 where $\mathcal{I}$ denotes the $k\times k$ unit matrix. For any $y\in\Omega$ and $\sigma\in\R^l$ with $0\in T(\sigma)$, we deduce from Lemmas~\ref{l.NC.product}, \ref{lemma.coderivative.product}, and \ref{prop.coderivative.calculus}
 that
 \begin{align}
 D^*\mathcal{P}\left((y,\sigma)|(y,0,0)\right)\left(\begin{array}{c}\alpha\\\beta\\\gamma\end{array}\right) & =\varDelta_{N_{(\Omega\times\{ 0\})}(y,0)}(-(\alpha,\beta))\times D^*T(\sigma|0)(\gamma),\nonumber\\
 ~ & =\varDelta(y,\alpha)\times D^*T(\sigma|0)(\gamma),\label{f.coder.P}
 \end{align}
 where $\varDelta:\R^k\times\R^k\rightrightarrows\R^k$ is defined by
\begin{align}\label{f.map.Delta.VS}
\varDelta(y,\alpha): =\left\{ 
\begin{array}{c l}
\{ 0\} & \textnormal{ if }0\in \alpha+ N_\Omega(y),\\
\varnothing & \textnormal{ if }0\notin \alpha+ N_\Omega(y).
\end{array}
\right.
\end{align}
 
\begin{proposition}\label{p.CS.mr}
In the setting of this section, let $(\bar x,\bar\sigma)$ be a solution to the coupled constraint system \eqref{f.CS}, and let $G,\mathcal{P}$ be taken from \eqref{f.CS.reform}. The following are equivalent:\vspace*{-0.05in}
\begin{enumerate}
\item[\bf(a)] $G+\mathcal{P}$ is metrically regular around $((\bar x,\bar\sigma),(0,0,0))$.

\item[\bf(b)] We have the implication
\begin{align}\label{f.mr.crit.CS}
\left.
\begin{array}{r l}
0 & \in \nabla_x\varPhi(\bar x,\bar\sigma)z+ N_\Omega(\bar x),\\
0 & \in \nabla_\sigma\varPhi(\bar x,\bar\sigma)z+D^*T(\bar\sigma|0)(\nu)\end{array}\right\}\Longrightarrow z=0,\;\nu=0.
\end{align}
 \item[\bf(c)] $T$ is metrically regular around $(\bar\sigma,0)$, and we have
\begin{align}\label{f.mr.crit.CS.II}
\left.
\begin{array}{r l}
0 & \in \nabla_x\varPhi(\bar x,\bar\sigma)z+ N_\Omega(\bar x),\\
0 & \in \nabla_\sigma\varPhi(\bar x,\bar\sigma)z+\rge(D^*T(\bar\sigma|0))
\end{array}\right\}\Longrightarrow z=0.
\end{align}
\end{enumerate}
\end{proposition}
{\bf Proof}. It follows from Lemma~\ref{lemma.coderivative.sum} together with \eqref{f.Part.der.G}--\eqref{f.map.Delta.VS} that
\begin{align*}
 ~ & D^*\left(G+\mathcal{P}\right)\left((\bar x,\bar\sigma)|(0,0,0)\right)\left(\begin{array}{c}\alpha\\\beta\\\gamma\end{array}\right) \\
 = & \left\{ \left(\begin{array}{c}
 -\alpha+\nabla_x\varPhi(\bar x,\bar\sigma)\beta\\
 \nabla_\sigma\varPhi(\bar x,\bar\sigma)\beta+\chi
 \end{array}\right)\left|\,
 \begin{array}{c}
 -\alpha\in N_{\Omega}(\bar x),\\
 \chi\in D^*T(\bar\sigma|\bar\nu)(\gamma) 
 \end{array}
 \right.\right\}.
 \end{align*}
 This verifies the equivalence between (a) and (b) by using Theorem~\ref{theorem.metric-regularity}.

 To justify (b)$\Rightarrow$(c), take $z=0$ in \eqref{f.mr.crit.CS} and get that the inclusion $0\in D^*T(\bar \sigma|0)(\nu)$ yields $\nu=0$. It follows from Theorem~\ref{theorem.metric-regularity} that
 $T$ is metrically regular around $(\bar \sigma,0)$. The condition \eqref{f.mr.crit.CS} trivially ensures the fulfillment of \eqref{f.mr.crit.CS.II}.

 To check finally implication (c)$\Rightarrow$(b), pick any $(z,\nu)$ such that the inclusions $0 \in \nabla_x\varPhi(\bar x,\bar\sigma)z+ N_\Omega(\bar x)$ and $0 \in \nabla_\sigma\varPhi(\bar x,\bar\sigma)z+D^*T(\bar\sigma|0)(\nu)$ hold simultaneously. The latter implies that $0 \in \nabla_\sigma\varPhi(\bar x,\bar\sigma)z+\rge(D^*T(\bar\sigma|0))$,  and so  \eqref{f.mr.crit.CS.II} yields $z=0$. As a consequence, we get $0\in D^*T(\bar\sigma|0)(\nu)$ because $(z,\nu)$ solves the inclusions above. Since $T$ is metrically regular around $(\bar\sigma,0)$, it follows from Theorem~\ref{theorem.metric-regularity} that $\nu=0$. \qed\vspace*{0.05in}
 
 Based on Theorem~\ref{theorem.metric-regularity}, the latter proposition can be used to draw conclusions about Lipschitzian stability properties of solutions to system \eqref{f.CS}. In what follows, we present sufficient conditions for metric 2-regularity and Gfrerer regularity. Invoking \eqref{f.Part.der.G} , the semiderivative of the gradient of $G$ at a point $(x,\sigma)\in\R^k\times\R^l$ for a direction $(w,\mu)\in\R^k\times\R^l$ is calculated by
 \begin{align}\label{f.secDer.G.VS}
(\nabla G)'((x,\sigma);(w,\mu)) =\left(\begin{array}{c c c}
 0 & (\nabla_x\varPhi)'((x,\sigma);(w,\mu)) & 0\\
 0 & (\nabla_\sigma\varPhi)'((x,\sigma);(w,\mu)) & 0
 \end{array}\right).
 \end{align}
 The next result concerns the case where $\Omega$ and $T$ are polyhedral.\vspace*{-0.05in}

\begin{proposition}\label{prop.m2r.SC.polyhedral}
 In the setting of this section, let $(\bar x,\bar\sigma)$ be a solution to
 \eqref{f.CS}, and let $G,\mathcal{P}$ be taken from \eqref{f.CS.reform}. Assume that $\Omega$ and $T$ are polyhedral. If $\nabla\varPhi$ is semidifferentiable at $(\bar x,\bar\sigma)$ for a unit direction $(w,\mu)\in\R^k\times\R^l$ with $DT(\bar\sigma|0)(\mu)\neq\varnothing$, then the following two statements are equivalent and ensure the metric 2-regularity of $G+\mathcal{P}$ around $((\bar x,\bar\sigma),(0,0,0))$ relative to $(w,\mu)$: \vspace*{-0.05in}
 \begin{enumerate}
\item[\bf(a)] We have the implication
 \begin{align*}
 \left.
 \begin{array}{r l}
 0 & \in \nabla_x\varPhi(\bar x,\bar\sigma)z+ N_\Omega(\bar x),\\
 0 & \in \nabla_\sigma\varPhi(\bar x,\bar\sigma)z+D^*T(\bar\sigma|0)(\nu),\\
 0 & \in (\nabla_x\varPhi)'((\bar x,\bar\sigma);(w,\mu))z+
 \nabla_x\varPhi(\bar x,\bar\sigma)\beta+N_\Omega(\bar x),\\
 0 & \in (\nabla_\sigma\varPhi)'((\bar x,\bar\sigma);(w,\mu))z+
 \nabla_\sigma\varPhi(\bar x,\bar\sigma)\beta+D^*T(\bar\sigma|0)(\gamma)
 \end{array}\right\}\;\Longrightarrow\;\begin{array}{c} z=0,\\\nu=0.\end{array}
\end{align*}\vspace*{-0.05in}
 \item[\bf(b)] $T$ is metrically regular around $(\bar \sigma,0)$, and it holds that\vspace*{-0.05in}
 \begin{align*}
 \left.
 \begin{array}{r l}
 0 & \in \nabla_x\varPhi(\bar x,\bar\sigma)z+ N_\Omega(\bar x),\\
 0 & \in \nabla_\sigma\varPhi(\bar x,\bar\sigma)z+\rge(D^*T(\bar\sigma|0)),\\
 0 & \in(\nabla_x\varPhi)'((\bar x,\bar\sigma);(w,\mu))z+
 \nabla_x\varPhi(\bar x,\bar\sigma)\beta+N_\Omega(\bar x),\\
 0 & \in(\nabla_\sigma\varPhi)'((\bar x,\bar\sigma);(w,\mu))z+
 \nabla_\sigma\varPhi(\bar x,\bar\sigma)\beta+\rge(D^*T(\bar\sigma|0))
 \end{array}\right\}\;\Longrightarrow\;z=0.
 \end{align*}
 \end{enumerate}\vspace*{-0.05in}
 In particular, assertions {\rm(a)} and {\rm(b)} are sufficient for the Gfrerer regularity of $G+\mathcal{P}$ around $((\bar x,\bar\sigma),(0,0,0))$ relative to $((w,\mu),\eta)$ for all $\eta\in\R^{k+p+q}$.
 \end{proposition}\vspace*{-0.1in}
 {\bf Proof}. To show first that (a) is sufficient for the claimed metric 2-regularity, observe that $\mathcal{P}$ is polyhedral due to the assumed polyhedrality of $\Omega$ and $T$. In particular, $DT(\bar\sigma|0)(\mu)\neq \varnothing$ yields $(w,\mu)\in\dom(D\mathcal{P}((\bar x,\bar\sigma)|(\bar y,\bar\nu))$. Thus it follows from Theorem~\ref{theorem.m2r.polyhedral} with $F:=G$ and $C:=\mathcal{P}$ that \eqref{f.m2r.map.polyhedral} is sufficient for the metric 2-regularity of $F+C=G+\mathcal{P}$ around $((\bar x,\bar\sigma),(0,0,0))$ relative to $(w,\mu)$. Thanks to \eqref{f.Part.der.G}--\eqref{f.map.Delta.VS} and \eqref{f.secDer.G.VS}, we see that \eqref{f.m2r.map.polyhedral} transforms into the condition in (a). 
 	
 Implication (a)$\Rightarrow$(b) can be checked in a standard way, and thus it remains to verify (b)$\Rightarrow$(a). To proceed, take a solution $(z,\nu)$ to the system in (a) and get $z=0$ by (b). The system in (a) is now reduced to
 \begin{align*}
 \begin{array}{l l}
 0 \in N_\Omega(\bar x), &  
 0 \in D^*T(\bar\sigma|0)(\nu),\\ 
 0 \in \nabla_x\varPhi(\bar x,\bar\sigma)\beta+N_\Omega(\bar x), &  
 0 \in \nabla_\sigma\varPhi(\bar x,\bar\sigma)\beta+D^*T(\bar\sigma|0)(\gamma),
 \end{array}
\end{align*}
 where the inclusions in the second line are independent of $\nu$, while
 the first inclusion holds trivially. By Theorem~\ref{theorem.metric-regularity}, the metric regularity
 of $T$ yields $0\in D^*T(\bar\sigma|0)(\nu)\Rightarrow\nu=0$. The final claim on Gfrerer regularity follows from here by
Theorem~\ref{theorem.m2r.gfrerer}.\qed\vspace*{0.05in}
 
Next we consider the unconstrained case of $\Omega=\R^k$ and deal with the modified mappings $\widetilde{G}:\R^l\times\R^l\rightarrow\R^{p+q}$ and $\widetilde{\mathcal{P}}:\R^k\times\R^l\rightrightarrows\R^{p+q}$ defined by
\begin{align}
\widetilde{G}(x,\sigma) & :=\left(\begin{array}{c}\varPhi(x,\sigma)\\ 0\end{array}\right),\quad
\widetilde{\mathcal{P}}(x,\sigma) & :=\left(\begin{array}{c} \{ 0\}\\ T(\sigma)\end{array}\right):=\left\{ \left.
\left(\begin{array}{c}0\\ b\end{array}\right)\right|\, b\in T(\sigma)
\right\}.\label{f.map.hP}
\end{align}
The unconstrained system in \eqref{f.CS} is equivalent to the generalized equation
 \begin{align*}
 0\in\widetilde{G}(x,\sigma)+\widetilde{\mathcal{P}}(x,\sigma),
\end{align*}
and thus the arguments above bring us to the expressions
\begin{align}\label{f.secDer.hG.VS}
\begin{array}{ll}
\nabla\widetilde{G}(x,\sigma)
=\left(\begin{array}{c c}
\nabla_x\varPhi(x,\sigma) & 0\\
\nabla_\sigma\varPhi(x,\sigma) & 0
\end{array}\right),\\
(\nabla \widetilde{G})'((x,\sigma);(w,\mu)) 
=\left(\begin{array}{c c}
(\nabla_x\varPhi)'((x,\sigma);(w,\mu)) & 0\\
(\nabla_\sigma\varPhi)'((x,\sigma);(w,\mu)) & 0
\end{array}\right),
\end{array}
\end{align}\vspace*{-0.2in}
\begin{align}\label{f.coder.hP}
D^*\widetilde{\mathcal{P}}\left((y,\sigma)|(0,0)\right)
\left(\begin{array}{c}\beta\\\gamma\end{array}\right)=
\{ 0\}\times D^*T(\sigma|0)(\gamma)
\end{align}
for suitable vectors $x,y,\sigma,w,\mu,\beta,\gamma$. This allows us to establish the metric 2-regularity and Gfrerer regularity in the case below, which is not covered by Proposition~\ref{prop.m2r.SC.polyhedral}.\vspace*{-0.05in}

\begin{proposition}\label{prop.m2r.SC.nonpolyhedral}
In the setting of this section, let $\Omega=\R^k$, let $(\bar x,\bar\sigma)$ be a solution to system \eqref{f.CS}, and let the mappings $\widetilde{G},\widetilde{\mathcal{P}}$ be taken from \eqref{f.map.hP}. Assume that $T$ is 
Lipschitz-like around $(\bar\sigma,0)$, and that we have
\begin{align}\label{f.SCNP.separate}
\textnormal{rge}(\nabla_\sigma\varPhi(\bar x,\bar\sigma))\bigcap\rge(D^*T(\bar\sigma|0))=\{ 0\}.
\end{align}
If $\nabla\varPhi$ is semidifferentiable at $(\bar x,\bar\sigma)$ for a unit direction $(w,\mu)\in\R^k\times\R^l$ with \linebreak
$DT(\bar\sigma)(\mu)\neq\varnothing$ and $T$ is metrically regular around $(\bar\sigma,0)$, 
then the following condition is sufficient for the metric 2-regularity of $\widetilde{G}+\widetilde{\mathcal{P}}$ around $((\bar x,\bar\sigma),(0,0))$ relative to $(w,\mu)$: the vector $z=0$ is the only $z\in\R^n$ satisfying the system
\begin{align}\label{f.SC.CS.nonpolyhedral}
\begin{array}{ll}
0=\nabla_x\varPhi(\bar x,\bar\sigma)z,\;
0\in\nabla_\sigma\varPhi(\bar x,\bar\sigma)z+\rge(D^*T(\bar\sigma|0)),\\
0\in(\nabla_x\varPhi)'((\bar x,\bar\sigma);(w,\mu))z +\nabla_x\varPhi(\bar x,\bar\sigma)\beta,\\
0\in(\nabla_\sigma\varPhi)'((\bar x,\bar\sigma);(w,\mu))z + 
\nabla_\sigma\varPhi(\bar x,\bar\sigma)\beta+\rge(D^*T(\bar\sigma|0)).
\end{array}
\end{align}
In particular, the latter is also sufficient for the Gfrerer regularity of $\widetilde{G}+\widetilde{\mathcal{P}}$ around $((\bar x,\bar\sigma),(0,0))$ relative to $((w,\mu),\eta)$ for any $\eta\in\R^{p+q}$.
\end{proposition}\vspace*{-0.1in}
{\bf Proof}. We intend to apply Theorem~\ref{theorem.m2r.productMap} with $\widetilde{G}$ and $\widetilde{\mathcal{P}}$. It follows from
\eqref{f.secDer.hG.VS} that $\nabla \widetilde{G}$ is semidifferentiable at $(\bar x,\bar\sigma)$ for $(w,\mu)$, and so Theorem~\ref{theorem.m2r.productMap}(a) holds. Conditions (b) and (c) of that theorem are guaranteed by the imposed assumption, and we get $\widetilde{\mathcal{P}}(x,\sigma)=R(x)\times T(\sigma)$ with $R\equiv \{ 0\}$ for all $(x,\sigma)$. Thus the first equality in (d) of Theorem~\ref{theorem.m2r.productMap} is satisfied. The second equality therein follows by \eqref{f.secDer.hG.VS} and \eqref{f.SCNP.separate}:
\begin{align*}
\textnormal{rge}(\nabla_\sigma\widetilde{G}(\bar x,\bar\sigma))\bigcap\rge(D^*T(\bar\sigma|0))=
\textnormal{rge}(\nabla_\sigma\varPhi(\bar x,\bar\sigma))\bigcap\rge(D^*T(\bar\sigma|0))=\{0\}.
\end{align*}
Finally, due to \eqref{f.coder.hP}, condition (e) of Theorem~\ref{theorem.m2r.productMap} transforms into
\begin{align*}
\left.
\begin{array}{r l}
0 & =\nabla_x\varPhi(\bar x,\bar\sigma)\theta,\\
0 & \in \nabla_\sigma\varPhi(\bar x,\bar\sigma)\theta+D^*T(\bar\sigma|0)(\zeta),\\
0 & \in (\nabla_x\varPhi)'((\bar x,\bar\sigma^;(w,\mu))\theta +
\nabla_x\varPhi(\bar x,\bar\sigma)\beta,\\
0 & \in (\nabla_\sigma\varPhi)'((\bar x,\bar\sigma^;(w,\mu))\theta + 
\nabla_\sigma\varPhi(\bar x,\bar\sigma)\beta+\rge(D^*T(\bar\sigma|0))
\end{array}\right\}\;\Rightarrow\;\begin{array}{c}
\theta=0,\\
\zeta=0.\end{array}
\end{align*}
By the metric regularity of $T$ around $(\bar \sigma,0)$, we see that the latter condition follows by the one in \eqref{f.SC.CS.nonpolyhedral}. Therefore, all the assumptions of Theorem~\ref{theorem.m2r.productMap} hold, and we arrive at the metric 2-regularity of $\widetilde{G}+\widetilde{\mathcal{P}}$ around $((\bar x,\bar\sigma),(0,0))$ relative to $(w,\mu)$. Theorem~\ref{theorem.m2r.gfrerer} ensures also the Gfrerer regularity relative to $((w,\mu),\eta)$ for all $\eta\in\R^{p+q}$.
\qed\vspace*{0.05in}

The metric regularity assumption on $T$ appears in both Propositions~\ref{prop.m2r.SC.polyhedral} and \ref{prop.m2r.SC.nonpolyhedral}. The price to pay for removing the polyhedrality assumption on $T$ is imposing the Lipschitz-like property in the unconstrained setting of \eqref{f.CS}.\vspace*{-0.2in}
 
\section{Applications to Variational Systems}\label{sec.App.VS}\vspace*{-0.1in}
 
 This section deals with (nonparametric) {\em variational systems} given in the form
 \begin{align}\label{f.VS}
 0\in f(x)+\mathcal{M}(x)N_{C_0}(g(x)),
 \end{align}
 where $f:\R^k\rightarrow\R^s$, $g:\R^k\rightarrow\R^q$, and  $\mathcal{M}:\R^k\rightarrow\R^{s\times q}$ are continuously differentiable, where $C_0\subset\R^q$ is closed, and where
\begin{align*}
\mathcal{M}(x)N_{C_0}(g(x)):=\bigcup_{\lambda\in N_{C_0}(g(x))}\mathcal{M}(x)\lambda.
 \end{align*}
A remarkable specialization of variational systems \eqref{f.VS} are {\em KKT systems}, where $\mathcal{M}(x):=\nabla g(x)$, and where $C_0$ is the convex polyhedral set $\R_{-}^{t}\times\{ 0\}$.

It is clear that \eqref{f.VS} is equivalently reformulated as
\begin{align}\label{f.VS.reform}
f(x)+\mathcal{M}(x)\lambda=0,\quad \zeta=g(x),\quad \lambda\in N_{C_0}(\zeta),
\end{align}
which can be written in turn in the constraint system form \eqref{f.CS} with $\Omega=\R^k$ and
\begin{align}
\varPhi(x,\sigma) & :=\varPhi(x,\lambda,\zeta):=\left(
\begin{array}{c}
f(x)+\mathcal{M}(x)\lambda\\
g(x)-\zeta
\end{array}\right),\label{f.map.Phi}\\
T(\sigma) & :=T(\lambda,\zeta):=-\lambda+N_{C_0}(\zeta).\label{f.map.T}
\end{align}
This allows us to rephrase the results of the previous section to guarantee the metric regularity, metric 2-regularity, and Gfrerer regularity of the summation mappings $G+\mathcal{P}$ and $\widetilde{G}+\widetilde{\mathcal{P}}$ defined above. It follows from Lemmas~\ref{lemma.coderivative.sum} and \ref{lemma.coderivative.product} with $\sigma=(\lambda,\zeta)$ that for all $\lambda\in N_{C_0}(\zeta)$ we have the equality
\begin{align}\label{f.cod.mapT}
D^*T(\sigma|0)(\alpha) & 
=\left\{ \left.\left(\begin{array}{c}-\alpha\\ \nu\end{array}\right)\right|\;\nu\in D^*N_{C_0}(\zeta|\lambda)(\alpha)\right\}.
\end{align}
The following lemma is an easy application of Theorem~\ref{theorem.metric-regularity}.\vspace*{-0.02in}

\begin{lemma}\label{l.mr.T} In the setting of this section, let $(\bar x,\bar\sigma):=(\bar x,\bar\lambda,\bar\zeta)$ be a solution to system \eqref{f.CS} with $\Omega=\R^k$ and with the mappings $\varPhi,T$ defined by \eqref{f.map.Phi} and \eqref{f.map.T}, respectively. Then $T$ is metrically regular around $(\bar\sigma=(\bar\lambda,\bar\zeta),0)$.
\end{lemma}\vspace*{-0.05in}

Now we apply the results of the preceding section to the case of \eqref{f.CS} with \eqref{f.map.Phi} and \eqref{f.map.T} to obtain efficient conditions for metric regularity, metric 2-regularity, and Gfrerer regularity via the given data. For $\varPhi$ in \eqref{f.map.Phi}, it is easy to see that
\begin{align}
\nabla_x\varPhi(x,\sigma) & =\nabla_x\varPhi(x,\lambda,\zeta)=\left(\begin{array}{c c}
\nabla f(x)+\sum\limits_{i=1}^q\lambda_i \nabla\mathcal{M}_i(x) & \;\;\nabla g(x)
\end{array}\right),\label{f.derx.map.Phi}\\
\nabla_\sigma\varPhi(x,\sigma) & =\nabla_{\lambda,\zeta}\varPhi(x,\lambda,\zeta)=\left(
\begin{array}{c c}
\mathcal{M}(x)^\top & 0\\
0 & -\mathcal{I}
\end{array}\right),\label{f.dersig.map.Phi}
\end{align}
where $\mathcal{I}$ is the $q\times q$ unit-matrix, and $\mathcal{M}_1(\cdot),\ldots,\mathcal{M}_q(\cdot)$ denoting the columns of $\mathcal{M}(\cdot)$. Then Proposition~\ref{p.CS.mr} induces the following one.\vspace*{-0.02in}

\begin{proposition}\label{prop.mr.VS}
In the setting of this section, let $(\bar x,\bar\sigma):=(\bar x,\bar\lambda,\bar\zeta)$ be a solution to system \eqref{f.CS} with $\Omega=\R^k$ and with $\varPhi,T$ defined in \eqref{f.map.Phi} and \eqref{f.map.T}. For $G,\mathcal{P}$ taken from \eqref{f.map.G} and \eqref{f.map.P}, the following statements are equivalent:\vspace*{-0.05in}
\begin{enumerate}
\item[\bf(a)] $G+\mathcal{P}$ is metrically regular around $((\bar x,\bar\sigma)=(\bar x,\bar\lambda,\bar\zeta),(0,0,0))$.
\item[\bf(b)] The vector $z=0$ is the only solution $z\in\R^s$ to the inclusion
\begin{align*}
0\in \left(\nabla f(\bar x)+\sum\limits_{i=1}^q\bar\lambda_i\nabla\mathcal{M}_i (\bar x)\right)z
+\nabla g(\bar x)D^*N_{C_0}(\bar\zeta|\bar\lambda)\left(\mathcal{M}(\bar x)^\top z\right).
\end{align*}
\end{enumerate}
\end{proposition}\vspace*{-0.05in}

The criterion for the metric regularity in the proposition uses the coderivative of the normal cone mapping to the set $C_0$, which is a particular case of the {\em second-order subdifferential/generalized Hessian} by Mordukhovich \cite{Mo92}. We refer the reader to the recent book \cite{Mo24} for the state-of-the-art results regarding this construction including its explicit computations for broad classes of extended-real-valued functions.\vspace*{0.03in}

We conclude this section with sufficient conditions for metric 2-regularity and Gfrerer regularity based on Proposition~\ref{prop.m2r.SC.polyhedral}. 
%and \ref{prop.m2r.SC.nonpolyhedral}. 
If $\nabla f,\nabla g,\nabla\mathcal{M}_1,\ldots,\nabla\mathcal{M}_q$ are semidifferentiable, it follows from \eqref{f.derx.map.Phi} and \eqref{f.dersig.map.Phi} with $\sigma=(\lambda,\zeta)$ that 
\begin{align*}
(\nabla_x\varPhi)'\left( (x,\sigma);(w,\alpha,v)\right)\mu & =
\left((\nabla f)'(x;w)+\sum_{i=1}^q\left(\lambda_i(\nabla\mathcal{M}_i)'(x;w)+\alpha_i\nabla \mathcal{M}_i(x)\right)\right)z \\ 
~ & \quad +(\nabla g)'(x;w)\chi,\\
(\nabla_\sigma\varPhi)'\left( (x,\sigma);(w,\alpha,v)\right)\mu & =\left(
\begin{array}{c}
\sum\limits_{j=1}^sz_j\nabla\widetilde{\mathcal{M}}_j(x)^\top w\\
0 
\end{array}\right)=
\left(
\begin{array}{c}
\sum\limits_{j=1}^sz_j(\widetilde{\mathcal{M}}_j)'(x) w\\
0 
\end{array}\right)
\end{align*}
for any $(\alpha,v)$ and $\mu=(z,\chi)$, where $\widetilde{\mathcal{M}}(\cdot):=\mathcal{M}(\cdot)^\top$, and $\widetilde{\mathcal{M}}_1(\cdot),\ldots,\widetilde{\mathcal{M}}_s(\cdot)$ are the columns of $\widetilde{\mathcal{M}}(\cdot)$. The next result follows from Lemma~\ref{l.mr.T} and Proposition~\ref{prop.m2r.SC.polyhedral}.
%\ref{prop.m2r.SC.nonpolyhedral}.
\vspace*{-0.05in}

\begin{proposition}\label{prop.m2r.VS.polyhedral}
In the setting of this section, let $(\bar x,\bar\sigma):=(\bar x,\bar\lambda,\bar\zeta)$ be a solution to system \eqref{f.CS} with $\Omega=\R^k$, and let the mappings $\varPhi,T$ be taken from \eqref{f.map.Phi}, \eqref{f.map.T}. Assume that $N_{C_0}$ is polyhedral $($which holds, in particular, when $C_0$ is polyhedral$)$ and consider such a direction $(w,v)\neq(0,0)$ with $DN_{C_0}(\bar\zeta|\bar\lambda)(v)\neq\varnothing$ that the mappings $\nabla f,\nabla g,\nabla\mathcal{M}_1,\ldots,\nabla\mathcal{M}_q$ are semidifferentiable at $\bar x$ for $w$. Then for $G,\mathcal{P}$ taken from \eqref{f.map.G} and \eqref{f.map.P}, the following conditions ensure the metric 2-regularity of $G+\mathcal{P}$ around $((\bar x,\bar\sigma)=(\bar x,\bar\lambda,\bar\zeta),(0,0,0))$ relative to $(w,\alpha,v)$ for all $\alpha\in\R^q$: \vspace*{-0.05in}
\begin{align*}
\left.
\begin{array}{r l}
0 & = \left(\nabla f(\bar x)+\sum\limits_{i=1}^q\bar\lambda_i\nabla\mathcal{M}_i(\bar x)\right)z
+\nabla g(\bar x)\chi,\\
\chi & \in D^*N_{C_0}(\bar\zeta|\bar\lambda)\left(\mathcal{M}(\bar x)^\top z\right),\\
0 & \in\left((\nabla f)'(\bar x;w)+\sum\limits_{i=1}^q\left(\lambda_i(\nabla\mathcal{M}_i)'(\bar x;w)+\alpha_i\nabla \mathcal{M}_i(\bar x)\right)\right)z \\
~ & \;\; +(\nabla g)'(\bar x;w)\chi+\left(\nabla f(\bar x)+\sum\limits_{i=1}^q\bar\lambda_i\nabla\mathcal{M}_i\bar (x)\right)\beta\\
~ & \;\;+\nabla g(\bar x)D^*N_{C_0}(\bar\zeta|\bar\lambda)\left(\mathcal{M}(\bar x)^\top \beta+\sum\limits_{j=1}^sz_j(\widetilde{\mathcal{M}}_j)'(\bar x)w\right)
\end{array}\right\}\,\Longrightarrow\,z=0,\chi=0.
\end{align*}
In particular, this condition yields the Gfrerer regularity of $G+\mathcal{P}$ around $((\bar x,\bar\sigma)=(\bar x,\bar\lambda,\bar\zeta),(0,0,0))$ relative to $((w,\alpha,v),\eta)$ for any $\alpha\in\R^q$ and any $\eta\in \R^{k+s+2q}$.
\end{proposition}\vspace*{-0.3in}

\section{Conclusions and Future Research}\label{sec:conc}\vspace*{-0.1in}

This paper develops a unified approach to the study of well-posedness properties in variational analysis and optimization, based on the notion of the least singular value function. This approach allows us to obtain novel results concerning broadly recognized concepts of metric regularity and Lipschitzian stability and also--which is the most important--to investigate newer notions of mixed-order regularity that require the usage of higher-order constructions of generalized differentiation. 

Since we are basically at the beginning of developing the new approach, many questions remain open for future research. They include constructive specifications of the general results obtained for particular classes of systems that appear overwhelmingly in variational analysis and optimization. Among them we mention parametric variational systems for which the standard metric regularity fails while the second-order regularity notions may be very useful. Another major direction of our future research is the application of the proposed regularity notions and results to higher-order numerical methods for optimization-related problems.\vspace*{-0.2in}

\bibliographystyle{spmpsci}

\end{document}